\documentclass[11pt]{article}

\usepackage[margin=1in]{geometry}
\usepackage{amsmath,amssymb,amsthm,mathtools}
\usepackage[T1]{fontenc}
\usepackage{lmodern}
\usepackage{microtype}
\usepackage{enumitem}
\usepackage[hidelinks,hypertexnames=false]{hyperref}
\usepackage{xcolor}
\definecolor{purple}{RGB}{112,48,160}
\usepackage{graphicx}
\usepackage{placeins}
\usepackage{float}
\graphicspath{{./}}
\usepackage{booktabs}
\numberwithin{equation}{section}

\makeatletter
\g@addto@macro\normalsize{%
  \setlength{\abovedisplayskip}{8pt plus 2pt minus 3pt}%
  \setlength{\belowdisplayskip}{8pt plus 2pt minus 3pt}%
  \setlength{\abovedisplayshortskip}{0pt plus 1pt}%
  \setlength{\belowdisplayshortskip}{5pt plus 1pt minus 2pt}%
}
\renewcommand\section{\@startsection{section}{1}{\z@}%
  {-2.7ex plus -.7ex minus -.2ex}%
  {1.5ex plus .2ex}%
  {\normalfont\Large\bfseries}}
\renewcommand\subsection{\@startsection{subsection}{2}{\z@}%
  {-2.2ex plus -.5ex minus -.2ex}%
  {1ex plus .2ex}%
  {\normalfont\large\bfseries}}
\makeatother

\newtheorem{theorem}{Theorem}
\newtheorem{lemma}[theorem]{Lemma}
\newtheorem{proposition}[theorem]{Proposition}

\theoremstyle{remark}

\newcommand{\E}{\mathbb E}

\newcommand{\R}{\mathbb R}
\newcommand{\T}{\mathbb T}
\newcommand{\cH}{\mathcal H}
\newcommand{\cL}{\mathcal L}
\newcommand{\norm}[1]{\left\lVert #1\right\rVert}

\newcommand{\dd}{\,\mathrm d}
\newcommand{\err}{\xi}
\newcommand{\avg}{\mathcal U}

\title{Strong Approximation of Stochastic Burgers Equations with Space--Time White Noise}
\author{Jianbo Cui{\thanks{{Department of Applied Mathematics, The Hong Kong Polytechnic University, Hung Hom, Hong Kong.}}}\; 
and Liying Sun{\thanks{{Academy for Multidisciplinary Studies, Capital Normal University, Beijing 100048, China.}}}}
\date{}

\begin{document}
\maketitle

\begin{abstract}
In this paper, we establish algebraic strong convergence rates for spectral Galerkin
and fully discrete approximations of the stochastic
Burgers equation driven by additive space--time white noise. For
deterministic initial data in $\cH^1$, the spectral approximation is shown to be 
convergent in the spatial supremum norm, uniformly in time and in every
finite probability moment, with order arbitrarily close to $1/2$.
Furthermore, we prove that an accelerated fully tamed exponential method possesses
 spatial and temporal strong convergence orders
arbitrarily close to $1/2$ and $1/4$  under a CFL condition. The key to our 
analysis is a polynomial stability estimate for the conservative Burgers
flux, which allows us to control the nonlinear error using polynomial
moment bounds for the exact and numerical solutions. To the best of our
knowledge, these are the first algebraic strong convergence rates on the
whole probability space for this equation with additive space--time white
noise.
\end{abstract}



\section{Introduction}
The randomly forced Burgers equation is a fundamental model for describing
turbulent transport, shock formation, and intermittency phenomena
\cite{BecKhanin2007}. Burgers-type models also arise in statistical physics
\cite{BecKhanin2007} and in stochastic adhesion descriptions of cosmological
structure formation \cite{MatarreseMohayaee2002}. For the space--time
white-noise problem considered here, closed-form solutions are generally
unavailable, making numerical approximation essential for investigating
solution trajectories and their statistics. However, constructing stable
approximations and establishing convergence rates are challenging because
the quadratic transport term is not globally Lipschitz and space--time
white noise limits the spatial and temporal regularity of the solution.

We consider the following stochastic Burgers equation
\begin{equation}
 \dd u(t,x)=\left[\nu\partial_{xx}u(t,x)
       +\partial_x\!\left(\frac{u(t,x)^2}{2}\right)\right]\dd t
       +\dd W(t,x),\qquad u(0,x)=u_0(x),\label{eq:burgers}
\end{equation}
on $\T_\ell=\R/(\ell\mathbb Z)$, where $\ell,\nu>0$ and $W$ is a
cylindrical Wiener process in $H=L^2(\T_\ell)$. The formal time derivative
of $W$ represents additive space--time white noise. For the corresponding
Dirichlet problem, global existence and uniqueness were established in
\cite{DaPratoDebusscheTemam1994} by subtracting the stochastic convolution
and analyzing the resulting random parabolic equation. Our goal is to
establish algebraic strong convergence rates in probability moments for
spatial and fully discrete approximations, uniformly on finite time
intervals.

In recent years, numerical approximation of stochastic Burgers equations
has received considerable attention. Almost-sure convergence rates for spatial finite difference
approximations driven by
white noise were proved in \cite{AlabertGyongy2006}. Convergence rates of
spectral Galerkin approximations for almost every noise realization were
established in \cite{BloemkerJentzen2013}. Such pathwise estimates take the
form
\begin{align}\label{path-estimate}
 \sup_{t\le T}\|u(t)-u_N(t)\|\le R(\omega)N^{-r},
\end{align}
where the finite random constant $R(\omega)$ depends on the noise
realization. For spatially correlated noise, pathwise convergence rates
for fully discrete approximations were obtained in
\cite{BloemkerKamraniHosseini2013}. For more general Burgers-type equations
driven by white noise, spatial convergence rates were proved up to
stopping times in \cite{HairerMatetski2016}.

These pathwise estimates describe the error for individual noise
realizations. Monte Carlo and multilevel Monte Carlo applications also
call for quantitative control in probability moments to assess accuracy
and computational cost \cite{Giles2015}. Such control does not follow
from the almost-sure finiteness of $R$ in \eqref{path-estimate} alone;
additional moment bounds are needed to obtain strong convergence rates. For trace-class noise, algebraic strong
convergence
rates for space--time discrete approximations were established for
additive noise in \cite{HutzenthalerJentzenLindnerPusnik2019} and
multiplicative noise in \cite{HutzenthalerLink2022}. Their proofs combine
nonlinear perturbation estimates with exponential moment bounds for the
exact and numerical solutions. Related strong-error techniques based on the exponential integrability also appear
in \cite{CuiHongLiu2017,HutzenthalerJentzen2020}. Recently, weak convergence rates
for spectral Galerkin approximations with additive
trace-class noise were obtained in \cite{BrehierCoxMillet2026} by
analyzing the associated Kolmogorov equation. 

The convergence analysis in the case of additive
space--time white noise is more delicate. Strong convergence of explicit
space--time approximations on the whole probability space was proved in
\cite{JentzenSalimovaWelti2019}, but the explicit algebraic rate was not established.
 Indeed, the exponential-moment arguments used for trace-class
noise can not be extended directly to this rougher setting. Recently, a fully discrete
tamed method for the related stochastic Burgers--Huxley equation with
space--time white noise was analyzed in \cite{WangCaoCao2026}, where one 
additional cubic reaction term is introduced to provide certain dissipation absent from
\eqref{eq:burgers}.

The main difficulties in this strong error analysis already appear in the spectral Galerkin
approximation
\[
 \dd u_N=\left[\nu\partial_{xx}u_N
       +P_N\partial_x\!\left(\frac{u_N^2}{2}\right)\right]\dd t
       +P_N\dd W,\qquad u_N(0)=P_Nu_0,
\]
where $P_N$ is the Fourier projection onto frequencies $|k|\le N$.
{The Fourier projections commute with spatial differentiation,
so the Galerkin approximation retains the conservative form of the nonlinear term.}
Preserving this structure is important under rough
forcing, since incompatible spatial discretizations can change the
limiting equation \cite{HairerVoss2011}. To isolate the nonlinear error,
write $S(t)=e^{\nu t\partial_{xx}}$ and
restore the omitted linear modes by setting
\[
 R_N(t)=(I-P_N)\left[S(t)u_0+\int_0^tS(t-s)\dd W(s)\right],
 \qquad \widehat u_N=u_N+R_N.
\]
Since $\widehat u_N$ and $u$ share the same initial value and driving
noise, their difference $e_N=u-\widehat u_N$ satisfies
\begin{align*}
 \partial_t e_N-\nu\partial_{xx}e_N
 &=\partial_x(\mathcal U_Ne_N+g_N),\qquad e_N(0)=0,\\
 \mathcal U_N&=\frac{u+\widehat u_N}{2},\qquad
 g_N=u_NR_N+\frac{R_N^2}{2}+(I-P_N)\!\left(\frac{u_N^2}{2}\right).
\end{align*}

The difficulty lies in controlling the nonlinear transport term. A
standard $L^2$ energy argument yields
\[
 \|\mathcal U_Ne_N\|_2\|\partial_xe_N\|_2
 \le \frac\nu2\|\partial_xe_N\|_2^2
     +C\|\mathcal U_N\|_2^4\|e_N\|_2^2.
\]
A direct Gronwall argument therefore gives a pathwise error estimate
with a random factor of the form
$\exp(C\int_0^T\|\mathcal U_N(t)\|_2^4\dd t)$.
Passing from this pathwise bound to a strong convergence rate requires
moment bounds for the exponential factor, which polynomial moment
estimates for the solutions alone do not provide. Even the quadratic
exponential integrability
of $\sup_{t\le T}\|u(t)\|_2^2$ established for Dirichlet boundary
conditions in \cite{DeVecchiJanakPriola2026}, including the white-noise
case, does not by itself control this fourth-power exponential.

We overcome this obstruction by using the conservative form of the
error equation to obtain an $L^1$ bound first:
\[
 \sup_{t\le T}\|e_N(t)\|_1
 \le C\bigl(1+\|\mathcal U_N\|_{L^2(0,T;H)}\bigr)
          \|g_N\|_{L^2(0,T;H)}.
\]
Combining this bound with one-dimensional interpolation closes the
$L^2$ energy estimate without introducing an exponential factor.
The resulting polynomial flux stability estimate,
Lemma~\ref{lem:stability} in Section~\ref{sec:stability}, gives
\[
 \sup_{t\le T}\|e_N(t)\|_2\le
 C\left[1+(1+\|\mathcal U_N\|_{L^2(0,T;H)})
                \|\mathcal U_N\|_{L^6(0,T;H)}^3\right]
       \|g_N\|_{L^2(0,T;H)}.
\]
The required polynomial moment bounds follow by subtracting a damped
stochastic convolution. For each moment order, a sufficiently large
damping parameter, fixed independently of $N$, reduces its covariance
enough to close the energy estimate using Gaussian exponential
integrability (see \eqref{eq:tag-4.5} and the proof of
Proposition~\ref{prop:continuous-inputs}). Combining these moment bounds
with spectral estimates for $g_N$ yields an $L^2$ strong error bound.
Heat-semigroup smoothing and the decomposition
$u-u_N=e_N+R_N$ then yield the supremum-norm rate. More precisely, for
deterministic $u_0\in \cH^1$, every $0<\gamma<1/2$, and
$2\le p<\infty$, Theorem~\ref{thm:spectral} gives
\[
 \left(\E\sup_{t\le T}\|u(t)-u_N(t)\|_\infty^p\right)^{1/p}
 \le C\lambda_{N+1}^{-\gamma/2},
\]
where $\lambda_N=1+(2\pi N/\ell)^2$. 

We then discretize the Galerkin equation in time by a tamed exponential
method that integrates the linear drift and samples the linear stochastic
convolution exactly \cite{JentzenKloedenWinkel2011}. The same noise
treatment appears in accelerated implicit approximations of the
stochastic Cahn--Hilliard equation
\cite{CuiHongSun2021}. Freezing and taming the nonlinear flux introduce
additional errors, which are controlled by the same polynomial stability
estimate. Let $Y_m^N$ denote the approximation at $t_m=m\tau$, with
$M\tau=T$. Under a CFL condition and the same assumptions on $u_0$,
$\gamma$, and $p$, Theorem~\ref{thm:tamed-linfty-grid} gives the coupled
bound
\[
 \left(\E\max_{0\le m\le M}\|u(t_m)-Y_m^N\|_\infty^p\right)^{1/p}
 \le C(\lambda_{N+1}^{-\gamma/2}+\tau^{\gamma/2}).
\]

In the rest of this paper, we first specify the schemes and state the main results in
Section~\ref{sec:main}. The proofs combine the Gaussian and analytic
estimates of Section~\ref{sec:tools} with the polynomial flux stability
bound developed in Section~\ref{sec:stability}. These tools yield the
spatial and fully discrete error estimates in
Sections~\ref{sec:spectral-proof} and~\ref{sec:tamed-exact}, respectively.
We conclude with numerical tests of the predicted convergence rates in
Section~\ref{sec:experiments}.

\section{Main results}\label{sec:main}
In this section, we introduce the basic notations and preliminaries, and then introduce the spatial and fully discrete approximations and state
our main results.

Set
\(
 H=L^2(\T_\ell),\quad L=-\partial_{xx},\quad
 S(t)=e^{-\nu tL},\quad \nu>0.
\)
The cylindrical Wiener process $W$ in \eqref{eq:burgers} is defined on a
filtered probability space
$(\Omega,\mathcal F,(\mathcal F_t)_{t\ge0},\mathbb P)$ satisfying the usual
conditions. With $F(r)=r^2/2$, the equation takes the abstract form
\[
 \dd u(t)+\nu L u(t)\dd t=\partial_xF(u(t))\dd t+\dd W(t),
 \qquad u(0)=u_0.
\]
We use the real periodic Fourier basis
\[
 e_0=\ell^{-1/2},\qquad
 e_{n,\cos}=\sqrt{2/\ell}\cos(2\pi nx/\ell),\qquad
 e_{n,\sin}=\sqrt{2/\ell}\sin(2\pi nx/\ell),
\]
with Laplace eigenvalues $\mu_n=(2\pi n/\ell)^2$, $n\ge0$.
The Galerkin projection $P_N$ retains the constant mode and complete
sine--cosine pairs up to frequency $N$. We write $H_N=P_NH$ and
$\lambda_N=1+\mu_N$. In sums over individual basis functions, $e_k$
denotes $e_0$, $e_{n,\cos}$, or $e_{n,\sin}$, with each positive
eigenvalue $\mu_n$ counted twice.
Throughout, $\cH^s=H^s(\T_\ell)$ denotes the inhomogeneous periodic
Sobolev space with norm
\[
 \|v\|_{\cH^s}^2
 =|\langle v,e_0\rangle|^2+
 \sum_{n\ge1}(1+\mu_n)^s
 (|\langle v,e_{n,\cos}\rangle|^2+|\langle v,e_{n,\sin}\rangle|^2).
\]
For negative $s$, this is understood by completion and duality.
We write $\|\cdot\|_q=\|\cdot\|_{L^q(\T_\ell)}$, so that
$\|\cdot\|_2$ is the norm of $H$, and use $\langle\cdot,\cdot\rangle$
for the $L^2$ inner product or the corresponding dual pairing.
Constants may depend on the fixed parameters $T,\nu,\ell,u_0$,
regularity exponents, and moment orders, but not on $N$ or $\tau$.
For the fully discrete estimates, constants may also depend on the fixed
$c_0$ in \eqref{eq:tag-6.0b}.
We assume $N\ge1$ and $0<\tau\le\min\{1,T\}$ throughout.
For a Banach space $X$, we use the abbreviations
$L^a_tX=L^a(0,T;X)$ and $C_tX=C([0,T];X)$, and set
$L^a_{t,x}=L^a((0,T)\times\T_\ell)$.

Setting $P_\infty=I$ and $u_\infty=u$ allows us to treat the exact
equation and its Galerkin approximations together:
\begin{equation}
 \dd u_N+\nu Lu_N\dd t=P_N\partial_xF(u_N)\dd t+P_N\dd W,
 \qquad u_N(0)=P_Nu_0,\qquad N\in\mathbb N\cup\{\infty\}.
 \label{eq:galerkin}
\end{equation}
Global mild solutions are obtained by subtracting the stochastic
convolution and studying the resulting random parabolic equation,
following the approach in \cite{DaPratoDebusscheTemam1994}.
Both drift terms have zero spatial integral, hence
\begin{equation}
 \langle u_N(t),e_0\rangle
 =\langle u_0,e_0\rangle+\langle W(t),e_0\rangle,
 \qquad N\in\mathbb N\cup\{\infty\}.
 \label{eq:shared-mean}
\end{equation}
We first derive norm bounds for smooth approximations by multiplying
the equations by the solution and integrating in space. Uniform bounds
then allow us to pass to the limit and justify these calculations for
the original solutions.

Our first main result gives an algebraic strong convergence rate for
the spectral Galerkin approximation.
\begin{theorem}\label{thm:spectral}
For deterministic $u_0\in\cH^1$, $2\le p<\infty$, and $0<\gamma<1/2$,
\begin{equation}
 \left(\E\sup_{t\le T}\|u(t)-u_N(t)\|_\infty^p\right)^{1/p}
 \le C\lambda_{N+1}^{-\gamma/2}. \label{eq:tag-5.8}
\end{equation}
\end{theorem}
Theorem~\ref{thm:spectral} gives an algebraic rate in probability moments,
complementing the pathwise rates in \cite{BloemkerJentzen2013} and the
strong convergence without an algebraic Burgers rate in
\cite{JentzenSalimovaWelti2019}. The key step is
Lemma~\ref{lem:stability}, which controls the effect of approximation
errors in the nonlinear term. Its bound depends only polynomially on
solution norms, so it can be used with their polynomial moment estimates.

We next discretize time on the grid $t_m=m\tau$, with $M\tau=T$.
Write
\[
 B_N(v)=P_N\partial_xF(v),\qquad S_N(t)=P_NS(t),\qquad
 \zeta_m^{N,\tau}=\int_{t_m}^{t_{m+1}}S_N(t_{m+1}-s)\dd W(s).
\]
Here $B_N$ is the projected Burgers drift and $\zeta_m^{N,\tau}$ is the
stochastic convolution over one time step. For $v\in H_N$, define
\begin{equation}
 \Theta_\tau(v)=
 \frac{1}{1+\tau^{1/4}\norm{B_N(v)}_{\cH^{-1}}},
 \qquad \widehat B_{N,\tau}(v)=\Theta_\tau(v)B_N(v).
 \label{eq:tag-6.0a}
\end{equation}
Combining this taming with the accelerated exponential approach in
\cite{JentzenKloeden2009,JentzenKloedenWinkel2011} gives
\begin{equation}
 \begin{gathered}
 Y_0^N=P_Nu_0,\\
 Y_{m+1}^N=S_N(\tau)Y_m^N
 +\int_0^\tau S_N(r)\widehat B_{N,\tau}(Y_m^N)\dd r+\zeta_m^{N,\tau}.
 \end{gathered}
 \label{eq:tag-6.0}
\end{equation}
The Fourier components of $\zeta_m^{N,\tau}$ are independent centered Gaussian
random variables. For a nonzero mode, the variance is
\[
 \frac{1-e^{-2\nu\mu_k\tau}}{2\nu\mu_k},
\]
and the constant-mode variance is $\tau$.

Our stability analysis uses the CFL condition
\begin{equation}
 \tau^{1/2}\lambda_N\le c_0,
 \label{eq:tag-6.0b}
\end{equation}
where $c_0>0$ is fixed independently of $N$ and $\tau$.
The taming factor bounds the drift in
$\cH^{-1}$, while this mesh condition also gives a bound in $H$; see
\eqref{eq:tag-6.2}--\eqref{eq:tag-6.3} in
Section~\ref{sec:tamed-exact}. These bounds are used to control the
moments of the numerical solution in Proposition~\ref{prop:tamed-moments}.
Our second main result gives an algebraic strong convergence rate for
\eqref{eq:tag-6.0}.

\begin{theorem}\label{thm:tamed-linfty-grid}
Let $u_0\in\cH^1$ be deterministic, $2\le p<\infty$, and
$0<\gamma<1/2$. Under \eqref{eq:tag-6.0b}, the scheme \eqref{eq:tag-6.0}
satisfies
\begin{equation}
 \left(\E\max_{0\le m\le M}\|u(t_m)-Y_m^N\|_\infty^p\right)^{1/p}
 \le C\left(\lambda_{N+1}^{-\gamma/2}+\tau^{\gamma/2}\right).
 \label{eq:main-full}
\end{equation}
\end{theorem}

Note that Theorems \ref{thm:spectral} and  \ref{thm:tamed-linfty-grid} also hold for the random initial data once the suitable polynomial order moment condition is imposed according to the proofs in Sections~\ref{sec:spectral-proof} and~\ref{sec:tamed-exact}.  We would also like to mention that implicit methods offer another way to stabilize the nonlinear drift,
but require solving nonlinear algebraic equations at each time step.
In an inexact implementation, the strong error analysis must also
account for errors from terminating Newton or fixed-point iterations,
while the nonlinear solves add to the computational cost.

\section{Preliminary estimates}\label{sec:tools}
We first recall the heat-semigroup and Sobolev estimates used to control
the nonlinear term, together with basic moment bounds for the stochastic
convolution. We then examine the effect of damping on its exponential
moments and derive the regularity and high-frequency tail estimates
needed in the error analysis.

\subsection{Semigroup and stochastic convolution estimates}
We collect standard heat-semigroup and Sobolev estimates, together
with moment and tail bounds for the stochastic convolution.

For $0<t\le T$, $1\le r\le q\le\infty$, and $a\ge b$, the periodic
heat kernel and analytic semigroup estimates (see, e.g.,  \cite{Pazy1983}, \cite[Chapter 6]{Ouhabaz2005}) give
\begin{align}
 \|\partial_x^jS(t)z\|_q
 &\le C t^{-j/2-\frac12(1/r-1/q)}\|z\|_r,
       \qquad j=0,1,\label{eq:heat-lpq}\\
 \|S(t)z\|_{\cH^a}
 &\le Ct^{-(a-b)/2}\|z\|_{\cH^b},\qquad
 \|\nu LS(t)z\|_{\cH^a}\le Ct^{-1-(a-b)/2}\|z\|_{\cH^b}.
 \label{eq:tag-4.8}
\end{align} We will also use the following periodic interpolation
and multiplication estimates (the latter follows by localization from
\cite[Theorem 5.1]{BehzadanHolst2021}):
\begin{align}
 \|v\|_4^8&\le C\|v\|_2^6\|v\|_{\cH^1}^2,
       \qquad \cH^{1/4}\hookrightarrow L^4,\label{eq:GN}\\
 \|F(v)\|_{\cH^\delta}&\le C\|v\|_{\cH^\beta}^2,
 \qquad \tfrac14<\beta<\tfrac12,\quad 0\le\delta<2\beta-\tfrac12.
 \label{eq:tag-2.10}
\end{align}

For a Hilbert space $X$, we use the time-fractional seminorm
\begin{equation}
 [v]_{W^{\theta,4}(I;X)}^4
 =\int_I\int_I\frac{\|v(t)-v(r)\|_X^4}{|t-r|^{1+4\theta}}\dd r\dd t.
 \label{eq:seminorm}
\end{equation}
For $1/4<\theta<1$, apply scalar fractional Morrey embedding
\cite[Theorem 8.2]{DiNezza2012} by duality, with the interval mean
removed using fractional Poincar\'e. Rescaling $I$ yields
\begin{equation}
 \|v(t)-v(r)\|_X\le C_\theta|t-r|^{\theta-1/4}
 [v]_{W^{\theta,4}(I;X)},\qquad t,r\in\overline I.
 \label{eq:morrey}
\end{equation} 
We also use the following standard spectral and $L^1$ heat estimates
\begin{align}
 \|(I-P_N)v\|_2&\le\lambda_{N+1}^{-s/2}\|v\|_{\cH^s},\quad s\ge0,
 \label{eq:tag-2.3}\\
 \|S(t)\partial_xh\|_{\cH^s}&\le Ct^{-3/4-s/2}\|h\|_1,
 \quad 0\le s<1/2.\label{eq:tag-2.4}
\end{align}
Here \eqref{eq:tag-2.4} follows from \eqref{eq:heat-lpq} and \eqref{eq:tag-4.8}.
For convenience, we write $Z(t)=\int_0^tS(t-r)\,\dd W(r)$.

\begin{lemma}[Stochastic convolution]\label{lem:convolution}
For $s<1/2$ and finite $p\ge1$,
\begin{align}
 \E\sup_{t\le T}\|Z(t)\|_{\cH^s}^p&<\infty,\label{eq:tag-2.5}\\
 \left(\E\sup_{t\le T}\|(I-P_N)Z(t)\|_2^p\right)^{1/p}
 &\le C\lambda_{N+1}^{-\gamma/2},\quad 0<\gamma<1/2.\label{eq:tag-2.6}
\end{align}
For $p\ge4$,
\begin{equation}
 \|(I-P_N)Z\|_{L^p(\Omega;L^4_tL^4_x)}
 \le C_p\lambda_{N+1}^{-1/4}.\label{eq:tag-2.7}
\end{equation}
\end{lemma}
\begin{proof}
The path-moment assertion \eqref{eq:tag-2.5} is the standard stochastic factorization
estimate \cite{DaPratoKwapienZabczyk1988,DaPratoZabczyk2014}. Then 
the spectral inequality \eqref{eq:tag-2.3} gives
\eqref{eq:tag-2.6} immediately.
To prove \eqref{eq:tag-2.7}, It\^o's isometry and the identity
$e_{n,\cos}(x)^2+e_{n,\sin}(x)^2=2/\ell$ give
\[
 \E|(I-P_N)Z(t,x)|^2
 =\frac{1}{\nu\ell}\sum_{n>N}
       \frac{1-e^{-2\nu\mu_nt}}{\mu_n}
 \le C\sum_{n>N}n^{-2}\le CN^{-1},
\]
uniformly in $t\in[0,T]$ and $x\in\T_\ell$.
For each fixed $(t,x)$, $(I-P_N)Z(t,x)$ is a centered Gaussian
random variable. Its $L^p(\Omega)$ norm is a constant depending only
on $p$ times its standard deviation. Consequently,
\[
 \bigl(\E|(I-P_N)Z(t,x)|^p\bigr)^{1/p}
 \le C_p\bigl(\E|(I-P_N)Z(t,x)|^2\bigr)^{1/2}
 \le C_pN^{-1/2}.
\]
It remains to integrate over time and space. Since $p/4\ge1$, the
triangle inequality for integrals in $L^{p/4}(\Omega)$
(Minkowski's integral inequality) yields
\begin{align*}
 &\left[\E\left(\int_0^T\int_{\T_\ell}
       |(I-P_N)Z(t,x)|^4\,\dd x\dd t\right)^{p/4}\right]^{1/p}\\
 &\qquad\le\left[\int_0^T\int_{\T_\ell}
       \bigl(\E|(I-P_N)Z(t,x)|^p\bigr)^{4/p}\,\dd x\dd t\right]^{1/4}
 \le C_p(T\ell)^{1/4}N^{-1/2}.\qedhere
\end{align*}
\end{proof}

\begin{lemma}[Deterministic linear tail]\label{lem:linear-tail}
For $u_0\in\cH^\gamma$, $\gamma\ge0$,
\begin{align}
 \sup_{t\le T}\|S(t)(I-P_N)u_0\|_2
 &\le\lambda_{N+1}^{-\gamma/2}\|u_0\|_{\cH^\gamma},\label{eq:tag-2.8}\\
 \|S(\cdot)(I-P_N)u_0\|_{L^4_tL^4_x}
 &\le C\lambda_{N+1}^{-\gamma/2-1/8}\|u_0\|_{\cH^\gamma}.
 \label{eq:tag-2.9}
\end{align}
\end{lemma}
\begin{proof}
The bound \eqref{eq:tag-2.8} follows from contractivity of $S(t)$ on $H$
and \eqref{eq:tag-2.3}. Since $(I-P_N)u_0$ contains only frequencies
$n>N$, \eqref{eq:GN} and spectral decay give
\[
 \|S(t)(I-P_N)u_0\|_4
 \le Ct^{-1/8}e^{-c\mu_{N+1}t}\|(I-P_N)u_0\|_2.
\]
Taking the fourth power and integrating in time, then applying
\eqref{eq:tag-2.3}, yields
\[
 \|S(\cdot)(I-P_N)u_0\|_{L^4_tL^4_x}
 \le C\mu_{N+1}^{-1/8}\|(I-P_N)u_0\|_2
 \le C\lambda_{N+1}^{-\gamma/2-1/8}\|u_0\|_{\cH^\gamma}.\qedhere
\]
\end{proof}
\subsection{Damped stochastic convolutions and exponential moments}
The solution energy estimates in Sections~\ref{sec:spectral-proof}
and~\ref{sec:tamed-exact} require exponential moments of the squared
supremum norm of the Gaussian part. We introduce damping to make the
admissible exponential parameter arbitrarily large. These Gaussian
estimates will yield polynomial moments of the nonlinear solutions and
approximations. 

For $N\in\mathbb N\cup\{\infty\}$, set $Z_{\alpha,N}(t)=\int_0^t e^{-(t-r)(\nu L+\alpha)}P_N\,\dd W(r)$,
$\alpha\ge1$,
which satisfies 
$$\dd Z_{\alpha,N}+\nu LZ_{\alpha,N}\dd t
 =-\alpha Z_{\alpha,N}\dd t+P_N\dd W.$$
  For each basis function retained by $P_N$, its modal variance is
\begin{equation}
 \E|\langle Z_{\alpha,N}(t),e_k\rangle|^2
 =\frac{1-e^{-2(\alpha+\nu\mu_k)t}}{2(\alpha+\nu\mu_k)}.
 \label{eq:tag-4.1a}
\end{equation}
The factorization arguments in
\cite[Sections~5.3--5.4]{DaPratoZabczyk2014} and Burkholder's inequality give that for every $1\le p<\infty$,
\begin{equation}
 \sup_{N,\alpha\ge1}\E\|Z_{\alpha,N}\|_{C_t\cH^\beta}^p<\infty,
 \qquad \beta<1/2.\label{eq:tag-4.1b}
\end{equation}
and 
\begin{equation}
 \sup_{N,\alpha,x}\|Z_{\alpha,N}(t,x)-Z_{\alpha,N}(r,x)\|_{L^p(\Omega)}
 \le C_p|t-r|^{1/4}.\label{eq:tag-4.6}
\end{equation}

We next make the dependence on $\alpha$ explicit. The resulting bounds
will allow us to choose $\alpha$ so that the exponential factors from
Gronwall's inequality have the moments needed to prove
\eqref{eq:continuous-energy-moments} and \eqref{eq:tag-6.9}.
Fix $0<s<1/2$ and choose $q_0=2/s>\max\{2,1/s\}$. Then by the  Sobolev embedding, we have 
\begin{equation}
 \|z\|_\infty\le C\|(I+L)^{s/2}z\|_{q_0}.\label{eq:tag-4.2}
\end{equation}
Notice that
\[
 \sum_{n\ge0}\frac{(1+\mu_n)^s}{\alpha+\nu\mu_n}
 \le C\alpha^{-1}\sum_{n\le\sqrt\alpha}(1+n^2)^s
       +C\sum_{n>\sqrt\alpha}n^{2s-2}
 \le C\alpha^{s-1/2}.
\]
Together with \eqref{eq:tag-4.1a}, this yields
\begin{equation}
 \sup_{N,t,x}\E|(I+L)^{s/2}Z_{\alpha,N}(t,x)|^2
 \le C\sum_{n\ge0}\frac{(1+\mu_n)^s}{\alpha+\nu\mu_n}
 \le C\alpha^{s-1/2}.\label{eq:tag-4.2a}
\end{equation}
Gaussian moment equivalence, Minkowski's inequality, and
\eqref{eq:tag-4.2} imply
\begin{equation}
 \sup_{N,t}\|\|Z_{\alpha,N}(t)\|_\infty\|_{L^q(\Omega)}
 \le C_s\sqrt q\,\alpha^{s/2-1/4},\qquad q\ge2.
 \label{eq:tag-4.3}
\end{equation} 
For $q<q_0$, we use H\"older's inequality.
Since $j!\ge(j/e)^j$ for $j\ge1$, expanding the exponential shows that, with
$c_s=(4eC_s^2)^{-1}$,
\begin{equation}
 \sup_{N,t}\E e^{\eta\|Z_{\alpha,N}(t)\|_\infty^2}\le2,
 \qquad 0\le\eta\le c_s\alpha^{1/2-s}.\label{eq:tag-4.4}
\end{equation}

Let $a>0$. Applying Jensen's inequality to the time average on $[0,T]$
and to the average over the $M$ grid points $t_0,\ldots,t_{M-1}$,
with $M\tau=T$, yields
\begin{equation}
 \max\left\{\E e^{a\int_0^T\|Z_{\alpha,N}(t)\|_\infty^2\dd t},
 \E e^{a\tau\sum_{m=0}^{M-1}\|Z_{\alpha,N}(t_m)\|_\infty^2}\right\}\le2
 \label{eq:tag-4.5}
\end{equation}
provided
\begin{equation}
 \alpha\ge\alpha_*(a):=\max\{1,(aT/c_s)^{2/(1-2s)}\}.
 \label{eq:tag-4.5a}
\end{equation}

Notice that by \eqref{eq:tag-4.1a}
\[
 \E\|\alpha Z_{\alpha,N}(t)\|_{\cH^{-1}}^2
 \le\frac{\alpha^2}{2}\sum_k\frac{(1+\mu_k)^{-1}}{\alpha+\nu\mu_k}
 \le C\alpha\sum_k(1+\mu_k)^{-1}\le C\alpha.
\]
Gaussian moments and Minkowski's inequality therefore give
\begin{equation}
 \|\alpha Z_{\alpha,N}\|_{L^p(\Omega;L^4_t\cH^{-1}_x)}
 \le C_{p,T}\alpha^{1/2}.\label{eq:tag-4.5b}
\end{equation} 
We shall choose $\alpha$ large enough to obtain the exponential moment bound
needed for the energy estimates \eqref{eq:continuous-energy-moments}
and \eqref{eq:tag-6.9}. 

\subsection{Deterministic convolutions and linear tails}\label{sec:deterministic-tails}

We first establish spatial regularity and time-increment estimates
for a forced heat equation. We then bound the Fourier truncation
error of the linear solution $S(t)u_0+Z(t)$. These estimates are used
in the strong error analysis in
Sections~\ref{sec:spectral-proof} and~\ref{sec:tamed-exact}.

For $u_0\in\cH^1$ and $f\in L^4(0,T;\cH^{-1})$, consider the heat equation
\[
 \partial_t\mathcal V_N+\nu L\mathcal V_N=P_Nf,
 \qquad \mathcal V_N(0)=P_Nu_0.
\]
The mild solution can be written as
\[
 \mathcal V_N(t)=S(t)P_Nu_0+
 \int_0^tS(t-s)P_Nf(s)\,\dd s.
\]

For $0\le\beta<1/2$, apply \eqref{eq:tag-4.8} with $a=\beta$, $b=-1$ and H\"older's inequality in time 
with exponents $4/3$ and $4$:
\begin{align*}
 \|\mathcal V_N(t)\|_{\cH^\beta}
 &\le\|u_0\|_{\cH^1}
       +C\int_0^t(t-s)^{-(1+\beta)/2}\|f(s)\|_{\cH^{-1}}\,\dd s\\
 &\le\|u_0\|_{\cH^1}
       +C\left(\int_0^t r^{-2(1+\beta)/3}\,\dd r\right)^{3/4}
           \|f\|_{L^4_t\cH^{-1}_x}.
\end{align*}
Taking the supremum over time gives
\begin{equation}
 \|\mathcal V_N\|_{C([0,T];\cH^\beta)}
 \le C_{T,\beta}\bigl(\|u_0\|_{\cH^1}+\|f\|_{L^4_t\cH^{-1}_x}\bigr),
 \qquad 0\le\beta<1/2.                              \label{eq:tag-4.9}
\end{equation}

We next estimate the increments $\mathcal V_N(t+h)-\mathcal V_N(t)$
in $\cH^{1/4}$. The embedding $\cH^{1/4}\hookrightarrow L^4$ in
\eqref{eq:GN} then gives the corresponding $L^4$ bounds.
For $0<h<T$ and $0<t<T-h$, subtracting the mild
formulas at $t+h$ and $t$ gives
\begin{align*}
 \mathcal V_N(t+h)-\mathcal V_N(t)
 &=S(t)(S(h)-I)P_Nu_0\\
 &\quad+\int_t^{t+h}S(t+h-s)P_Nf(s)\,\dd s\\
 &\quad+\int_0^t[S(t+h-s)-S(t-s)]P_Nf(s)\,\dd s.
\end{align*}
Since $1-e^{-\nu\mu_kh}\le(\nu\mu_kh)^{3/8}$, the Fourier expansion gives
\[
 \|S(\cdot)(S(h)-I)P_Nu_0\|_{L^4(0,T-h;\cH^{1/4})}
 \le C_T h^{3/8}\|u_0\|_{\cH^1}.
\]
By \eqref{eq:tag-4.8}, the triangle inequality and the identity
$S(r+h)-S(r)=-\int_0^h\nu LS(r+a)\,\dd a$ give, for $r>0$,
\begin{align*}
 \|S(r)\|_{\cL(\cH^{-1},\cH^{1/4})}
 &\le Cr^{-5/8},\\
 \|S(r+h)-S(r)\|_{\cL(\cH^{-1},\cH^{1/4})}
 &\le C\min\{r^{-5/8},hr^{-13/8}\}.
\end{align*}
Extending $f$ by zero outside $(0,T)$, we apply Young's inequality
in time to the integrals over $(t,t+h)$ and $(0,t)$ to obtain
\begin{align*}
 &\left\|\int_t^{t+h}S(t+h-s)P_Nf(s)\,\dd s
       \right\|_{L^4(0,T-h;\cH^{1/4})}\\
 &\qquad\le C\left(\int_0^h r^{-5/8}\,\dd r\right)
       \|f\|_{L^4(0,T;\cH^{-1})}
 \le Ch^{3/8}\|f\|_{L^4(0,T;\cH^{-1})},\\
 &\left\|\int_0^t[S(t+h-s)-S(t-s)]P_Nf(s)\,\dd s
       \right\|_{L^4(0,T-h;\cH^{1/4})}\\
 &\qquad\le C\left(\int_0^h r^{-5/8}\,\dd r
       +h\int_h^T r^{-13/8}\,\dd r\right)
       \|f\|_{L^4(0,T;\cH^{-1})}
 \le Ch^{3/8}\|f\|_{L^4(0,T;\cH^{-1})}.
\end{align*}
The triangle inequality in the increment decomposition now gives
\begin{equation}
 \|\mathcal V_N(\cdot+h)-\mathcal V_N\|_{L^4(0,T-h;\cH^{1/4})}
 \le C_T h^{3/8}
       \bigl(\|u_0\|_{\cH^1}+\|f\|_{L^4_t\cH^{-1}_x}\bigr).
 \label{eq:tag-4.10}
\end{equation}
 By symmetry and the change of variables $h=t-r$, using \eqref{eq:tag-4.10},
for $0<\vartheta<3/8$ we obtain
\begin{align*}
 [\mathcal V_N]_{W^{\vartheta,4}(0,T;\cH^{1/4})}^4
 &=2\int_0^T h^{-1-4\vartheta}
       \|\mathcal V_N(\cdot+h)-\mathcal V_N\|_{L^4(0,T-h;\cH^{1/4})}^4\,\dd h\\
 &\le C\bigl(\|u_0\|_{\cH^1}+\|f\|_{L^4_t\cH^{-1}_x}\bigr)^4
       \int_0^T h^{3/2-1-4\vartheta}\,\dd h<\infty.
\end{align*}
Combining this with \eqref{eq:tag-4.9} at $\beta=1/4$ yields
\begin{equation}
 \|\mathcal V_N\|_{W^{\vartheta,4}(0,T;\cH^{1/4})}
 \le C_{T,\vartheta}
       \bigl(\|u_0\|_{\cH^1}+\|f\|_{L^4_t\cH^{-1}_x}\bigr),
 \qquad 0<\vartheta<3/8.
 \label{eq:tag-4.11}
\end{equation}


We next estimate the error from truncating the Fourier expansion of
the linear solution $S(t)u_0+Z(t)$.

\begin{lemma}[Linear high-frequency tail estimates]
\label{lem:linfty-tail}
Let $u_0\in\cH^1$, $1\le p<\infty$, and $0<\gamma<1/2$.  For $N\in\mathbb N$ and $K\in\mathbb N\cup\{\infty\}$ with $K>N$, put
\[
 R_{N,K}(t)=S(t)(P_K-P_N)u_0+(P_K-P_N)Z(t).
\]
For $K=\infty$, write $R_N=R_{N,\infty}$. Then, uniformly in $K>N$,
\begin{equation}
 \|R_{N,K}\|_{L^p(\Omega;C_tH)}
 +\|R_{N,K}\|_{L^p(\Omega;L^4_{t,x})}
 +\|R_{N,K}\|_{L^p(\Omega;C_tL^\infty_x)}
 \le C_{p,\gamma}\lambda_{N+1}^{-\gamma/2}.
 \label{eq:tag-7.1}
\end{equation}
\end{lemma}

\begin{proof}
Since $1-\gamma>1/2$, the embedding
$\cH^{1-\gamma}\hookrightarrow L^\infty$ and the spectral estimate
\eqref{eq:tag-2.3} give
\[
 \sup_{t\le T}\|S(t)(P_K-P_N)u_0\|_\infty
 \le C_\gamma\|(P_K-P_N)u_0\|_{\cH^{1-\gamma}}
 \le C_\gamma\lambda_{N+1}^{-\gamma/2}\|u_0\|_{\cH^1}.
\]
For the stochastic convolution, we use the factorization method \cite{DaPratoKwapienZabczyk1988}; see also
\cite[Lemma~4.3]{BloemkerJentzen2013} for the corresponding estimate
with sine modes. Choose $s,\theta>0$ such that
$s+2\theta<1/2-\gamma$, and then choose
$q>\max\{p,2,1/s,1/\theta\}$. Set
\[
 Y(t)=\int_0^t(t-r)^{-\theta}S(t-r)(P_K-P_N)\,\dd W(r).
\]
It\^o's isometry and the boundedness of the Fourier basis functions give
\begin{align*}
 \sup_{t\le T,\,x\in\T_\ell}\E|(I+L)^{s/2}Y(t,x)|^2
 &\le C\sum_{n>N}(1+\mu_n)^s
          \int_0^T r^{-2\theta}e^{-2\nu\mu_nr}\,\dd r\\
 &\le C\sum_{n>N}n^{2s+4\theta-2}
 \le CN^{2s+4\theta-1}.
\end{align*}
The Gaussian moment bound at each $(t,x)$, integrated over time and
space, therefore yields
\[
 \left(\E\int_0^T\|(I+L)^{s/2}Y(t)\|_q^q\,\dd t\right)^{1/q}
 \le C_qN^{s+2\theta-1/2}.
\]
The factorization identity reads
\[
 (P_K-P_N)Z(t)
 =\frac{\sin(\pi\theta)}{\pi}
   \int_0^t(t-r)^{\theta-1}S(t-r)Y(r)\,\dd r.
\]
Since $s>1/q$, Sobolev embedding and $L^q$ contractivity of $S(t)$
bound the integrand in the spatial supremum norm. Since $\theta>1/q$,
H\"older's inequality in time then gives
\begin{align*}
 \sup_{t\le T}\|(P_K-P_N)Z(t)\|_\infty
 &\le C\sup_{t\le T}\int_0^t(t-r)^{\theta-1}
          \|(I+L)^{s/2}Y(r)\|_q\,\dd r\\
 &\le C\left(\int_0^T\|(I+L)^{s/2}Y(r)\|_q^q\,\dd r\right)^{1/q}.
\end{align*}
Taking $L^p(\Omega)$ norms and using $p<q$ and the Gaussian moment
bound above, we obtain
\[
 \left(\E\sup_{t\le T}\|(P_K-P_N)Z(t)\|_\infty^p\right)^{1/p}
 \le CN^{s+2\theta-1/2}
 \le C\lambda_{N+1}^{-\gamma/2}.
\]
Combining this with the initial-data estimate gives the bound in
$L^p(\Omega;C_tL^\infty_x)$. The other two bounds follow from
\[
 \|R_{N,K}\|_{C_tH}+\|R_{N,K}\|_{L^4_{t,x}}
 \le\bigl(\ell^{1/2}+(T\ell)^{1/4}\bigr)
       \|R_{N,K}\|_{C_tL^\infty_x}.
\]
This proves \eqref{eq:tag-7.1}.
\end{proof}

\section{Polynomial flux stability}\label{sec:stability}

We now prove the polynomial stability estimate used in
Sections~\ref{sec:spectral-proof} and~\ref{sec:tamed-exact}.

\begin{lemma}[Polynomial flux stability]\label{lem:stability}
Fix $\nu>0$. Let $\err\in C([0,T];H)\cap L^2(0,T;\cH^1)$ satisfy
\begin{equation}
 \partial_t\err-\nu\partial_{xx}\err
 =\sigma\partial_x(\avg\err+g),\qquad \err(0)=0,                 \label{eq:tag-3.1}
\end{equation}
where \(\sigma\in\R\),
\(
 \avg\in L^6(0,T;H),
 g\in L^2(0,T;H).
\)
Assume that $\err$ is spatially periodic and $\int_{\T_\ell}\err(t,x)\dd x=0$ for almost every $t$.
Then

\begin{equation}
 \sup_{t\le T}\norm{\err(t)}_2
 \le C_{\nu,\ell,\sigma}
 \left[1+\bigl(1+\norm \avg_{L^2_tH}\bigr)
              \norm \avg_{L^6_tH}^3\right]
 \norm g_{L^2_tH}.                                 \label{eq:tag-3.2}
\end{equation}
\end{lemma}

\begin{proof}
{
It suffices to prove the estimate for $\sigma=1$. We only carry out the argument for smooth $\avg,g,\err$ since one can justify the
passage to weak solutions by standard procedures.  Put
\[
 \phi_\varepsilon(r)=\sqrt{r^2+\varepsilon^2}-\varepsilon,
 \qquad \varepsilon>0.
\]
A direct calculation gives
\[
 \phi_\varepsilon'(r)=\frac{r}{\sqrt{r^2+\varepsilon^2}},\qquad
 \phi_\varepsilon''(r)=
   \frac{\varepsilon^2}{(r^2+\varepsilon^2)^{3/2}},
\]
so that
\[
 0\le\phi_\varepsilon''(r)\le\varepsilon^{-1},\qquad
 r^2\phi_\varepsilon''(r)\le C\varepsilon,
 \qquad |r|\le\phi_\varepsilon(r)+\varepsilon.
\]
Multiplying \eqref{eq:tag-3.1} by $\phi_\varepsilon'(\err)$ and integrating by parts gives
\begin{align*}
 \frac{\dd}{\dd t}\int_{\T_\ell}\phi_\varepsilon(\err)\dd x
 +\nu\int_{\T_\ell}\phi_\varepsilon''(\err)|\err_x|^2\dd x
 =-\int_{\T_\ell}(\avg\err+g)\phi_\varepsilon''(\err)\err_x\dd x.
\end{align*} Then weighted
Young's inequality gives
\begin{align*}
 |(\avg\err+g)\phi_\varepsilon''(\err)\err_x|
 &\le \frac\nu2\phi_\varepsilon''(\err)|\err_x|^2
   +C_\nu(\avg^2\err^2+g^2)\phi_\varepsilon''(\err),
\end{align*}
and therefore
\begin{align*}
 \frac{\dd}{\dd t}\int_{\T_\ell}\phi_\varepsilon(\err)\dd x
 +\frac\nu2\int_{\T_\ell}\phi_\varepsilon''(\err)|\err_x|^2\dd x
 \le C_\nu\varepsilon\norm \avg_2^2
      +C_\nu\varepsilon^{-1}\norm g_2^2.
\end{align*}
Since $\err(0)=0$, integration in time and
$\norm \err_1\le\int_{\T_\ell}\phi_\varepsilon(\err)\dd x+\ell\varepsilon$ give
\[
 \sup_{t\le T}\norm{\err(t)}_1
 \le C_{\nu,\ell,\sigma}\left[
 \varepsilon\bigl(1+\norm \avg_{L^2_tH}^2\bigr)
 +\varepsilon^{-1}\norm g_{L^2_tH}^2\right].
\]
If $g\ne0$, choose
$\varepsilon=\norm g_{L^2_tH}/(1+\norm \avg_{L^2_tH})$; if $g=0$,
let $\varepsilon\downarrow0$.  In both cases
\begin{equation}
 \sup_{t\le T}\norm{\err(t)}_1
 \le C_{\nu,\ell,\sigma}\bigl(1+\norm \avg_{L^2_tH}\bigr)
                  \norm g_{L^2_tH}.                 \label{eq:tag-3.3}
\end{equation}

Testing \eqref{eq:tag-3.1} by $\err$ yields
\begin{align*}
 \frac12\frac{\dd}{\dd t}\norm \err_2^2+\nu\norm{\err_x}_2^2
 &=-\int_{\T_\ell}(\avg\err+g)\err_x\dd x\\
 &\le \norm \avg_2\norm \err_\infty\norm{\err_x}_2
      +\norm g_2\norm{\err_x}_2.
\end{align*}
Using the Nash--Gagliardo--Nirenberg inequality 
\begin{equation}
 \|\err\|_\infty\le C_\ell\|\err\|_1^{1/3}\|\err_x\|_2^{2/3},
 \label{eq:tag-3.4}
\end{equation} we have 
\[
 \norm \avg_2\norm \err_\infty\norm{\err_x}_2
 \le C_\ell\norm \avg_2\bigl(\sup_{r\le T}\|\err(r)\|_1\bigr)^{1/3}\norm{\err_x}_2^{5/3}.
\]
Together with Young's inequality, this gives
\begin{align*}
 \frac12\frac{\dd}{\dd t}\norm\err_2^2+\nu\norm{\err_x}_2^2
 &\le C_\ell\norm\avg_2
       \bigl(\sup_{r\le T}\|\err(r)\|_1\bigr)^{1/3}
       \norm{\err_x}_2^{5/3}
       +\norm g_2\norm{\err_x}_2\\
 &\le \frac\nu2\norm{\err_x}_2^2
       +C_{\nu,\ell,\sigma}\norm\avg_2^6
       \bigl(\sup_{r\le T}\|\err(r)\|_1\bigr)^2
       +C_\nu\norm g_2^2.
\end{align*}
Hence, by Gronwall's inequality, 
\[
 \sup_{t\le T}\norm{\err(t)}_2^2
 \le C_{\nu,\ell,\sigma}
 \left(\bigl(\sup_{r\le T}\|\err(r)\|_1\bigr)^2\norm \avg_{L^6_tH}^6+\norm g_{L^2_tH}^2\right).
\]
Taking square roots and inserting \eqref{eq:tag-3.3} proves \eqref{eq:tag-3.2}.
}
\end{proof}

\section{Solution moments and spectral convergence}\label{sec:spectral-proof}
We first establish moment bounds for the exact solution and its spectral
Galerkin approximations, uniformly in $N$. These bounds, together with
the linear tail estimate in Lemma~\ref{lem:linfty-tail}, allow us to apply
Lemma~\ref{lem:stability} and obtain a strong error bound in $H$.
Applying the heat-semigroup estimate \eqref{eq:heat-lpq} to the mild error
equation then yields the supremum-norm convergence rate in
Theorem~\ref{thm:spectral}.

\begin{proposition}[Solution moments]\label{prop:continuous-inputs}
Let $u_0\in\cH^1$. For every finite $p\ge1$ and $\beta<1/2$,
\begin{equation}
 \sup_{N\in\mathbb N\cup\{\infty\}}
 \E\left[\sup_{t\le T}\|u_N(t)\|_{\cH^\beta}^p
              +\|u_N\|_{L^8_tL^4_x}^p\right]<\infty.
 \label{eq:continuous-moments}
\end{equation}
\end{proposition}
\begin{proof}
Set $v_N=u_N-Z_{\alpha,N}$. Then
\[
 v_N'+\nu Lv_N=P_N\partial_xF(v_N+Z_{\alpha,N})+\alpha Z_{\alpha,N},
 \qquad v_N(0)=P_Nu_0.
\]
Testing against $v_N$ removes $P_N$. The identity
$\int_{\T_\ell}v_N^2\partial_xv_N\dd x=0$ and Young's inequality give
\begin{align}
 \frac{\dd}{\dd t}\|v_N\|_2^2+\nu\|\partial_xv_N\|_2^2
 &\le C_G(1+\|Z_{\alpha,N}\|_\infty^2)\|v_N\|_2^2\notag\\
 &\quad+C\|Z_{\alpha,N}\|_\infty^4+C\alpha^2\|Z_{\alpha,N}\|_2^2.
 \label{eq:continuous-energy}
\end{align}
Here $C_G\ge1$ depends on $\nu,\ell,T$, but not on $N,\alpha$.
Gronwall's inequality and integration of \eqref{eq:continuous-energy} yield
\begin{align*}
 &\sup_{t\le T}\|v_N(t)\|_2^2+\int_0^T\|v_N(t)\|_{\cH^1}^2\dd t\\
 &\quad\le C\left(1+\|u_0\|_2^2+
 \int_0^T[\|Z_{\alpha,N}\|_\infty^4+\alpha^2\|Z_{\alpha,N}\|_2^2]\dd t\right)
 \exp\!\left(C_GT+C_G\int_0^T\|Z_{\alpha,N}\|_\infty^2\dd t\right),
\end{align*}
after enlarging $C_G$ if needed. For a target moment $p$, choose
$\alpha\ge\alpha_*(16pC_G)$ from \eqref{eq:tag-4.5a}.
H\"older's inequality, \eqref{eq:tag-4.3}, \eqref{eq:tag-4.1b}, and
\eqref{eq:tag-4.5} give
\begin{equation}
 \|v_N\|_{L^{2p}(\Omega;C_tH)}
 +\|v_N\|_{L^{2p}(\Omega;L^2_t\cH^1)}\le C_p(1+\alpha).
 \label{eq:continuous-energy-moments}
\end{equation}
By \eqref{eq:GN},
\[
 \|v_N\|_{L^8_tL^4_x}^8
 \le C\sup_t\|v_N(t)\|_2^6\int_0^T\|v_N\|_{\cH^1}^2\dd t.
\]
Taking $L^{2p}(\Omega)$ norms, using \eqref{eq:continuous-energy-moments},
and adding $Z_{\alpha,N}$ with \eqref{eq:tag-4.1b} gives
$\|u_N\|_{L^{2p}(\Omega;L^8_tL^4_x)}\le C_p(1+\alpha)$. By $\|P_N\partial_xF(u_N)\|_{\cH^{-1}}\le C\|u_N\|_4^2$ and \eqref{eq:tag-4.5b},
\begin{equation}
 \|P_N\partial_xF(u_N)+\alpha Z_{\alpha,N}\|_{L^p(\Omega;L^4_t\cH^{-1})}
 \le C_p\bigl((1+\alpha)^2+\alpha^{1/2}\bigr).
 \label{eq:continuous-forcing}
\end{equation}
Apply \eqref{eq:tag-4.9} with $\mathcal V_N=v_N$ and
$f=P_N\partial_xF(u_N)+\alpha Z_{\alpha,N}$. For $0\le\beta<1/2$,
taking $L^p(\Omega)$ norms and using \eqref{eq:continuous-forcing} gives
\begin{align*}
 \|v_N\|_{L^p(\Omega;C_t\cH^\beta)}
 &\le C_{T,\beta}\left(\|u_0\|_{\cH^1}
 +\|P_N\partial_xF(u_N)+\alpha Z_{\alpha,N}\|_{L^p(\Omega;L^4_t\cH^{-1})}\right)\\
 &\le C_{p,\beta}(1+\alpha)^2.
\end{align*}
Since $u_N=v_N+Z_{\alpha,N}$, \eqref{eq:tag-4.1b} gives the corresponding
$C_t\cH^\beta$ moment bound for $u_N$. We fix
$\alpha=\alpha_*(16pC_G)$, independently of $N$. Together with the
$L^8_tL^4_x$ bound obtained above, this proves
\eqref{eq:continuous-moments} for $0\le\beta<1/2$.
\end{proof}

Now we are in a position to prove Theorem~\ref{thm:spectral}. 

\begin{proof}[\textbf{Proof of Theorem~\ref{thm:spectral}}]
Set $\widehat u_N=u_N+R_N$, $\err_N=u-\widehat u_N$, and
$\avg_N=(u+\widehat u_N)/2$. Subtracting the equations for $u$ and
$\widehat u_N$ gives
\begin{align*}
 \partial_t\err_N+\nu L\err_N
 &=\partial_x(\avg_N\err_N+g_N),\qquad \err_N(0)=0,\\
 g_N&=u_NR_N+\frac{R_N^2}{2}+(I-P_N)F(u_N).
\end{align*}
The definition of $R_N$ and \eqref{eq:shared-mean}  also show
that $\err_N$ has zero spatial mean, as required by
Lemma~\ref{lem:stability}.

To estimate the term $(I-P_N)F(u_N)$ in $g_N$, fix $0<\gamma<1/2$ and choose
\begin{equation}
 \gamma<\delta<\frac12,
 \qquad \frac{\delta+1/2}{2}<\beta<\frac12.
 \label{eq:tag-5.7}
\end{equation}
Then $\delta<2\beta-1/2$, so the spectral estimate
\eqref{eq:tag-2.3} and the multiplication estimate
\eqref{eq:tag-2.10} give
\[
 \|(I-P_N)F(u_N(t))\|_2
 \le\lambda_{N+1}^{-\delta/2}\|F(u_N(t))\|_{\cH^\delta}
 \le C\lambda_{N+1}^{-\delta/2}\|u_N(t)\|_{\cH^\beta}^2.
\]
Together with the expression for $g_N$, this yields
\[
 \|g_N(t)\|_2\le
 \|u_N(t)\|_2\|R_N(t)\|_\infty+C\|R_N(t)\|_\infty^2
 +C\lambda_{N+1}^{-\delta/2}\|u_N(t)\|_{\cH^\beta}^2.
\]
Taking the supremum over time and applying H\"older's inequality, we obtain from \eqref{eq:continuous-moments} and
\eqref{eq:tag-7.1}
\begin{equation}
 \begin{aligned}
 \|g_N\|_{L^p(\Omega;L^\infty_tH)}
 &\le \|u_N\|_{L^{2p}(\Omega;C_tH)}
       \|R_N\|_{L^{2p}(\Omega;C_tL^\infty_x)}\\
 &\quad+C\|R_N\|_{L^{2p}(\Omega;C_tL^\infty_x)}^2
       +C\lambda_{N+1}^{-\delta/2}
       \|u_N\|_{L^{2p}(\Omega;C_t\cH^\beta)}^2\\
 &\le C\bigl(\lambda_{N+1}^{-\gamma/2}
       +\lambda_{N+1}^{-\gamma}
       +\lambda_{N+1}^{-\delta/2}\bigr)
 \le C\lambda_{N+1}^{-\gamma/2}.
 \end{aligned}
 \label{eq:spectral-flux}
\end{equation}

By \eqref{eq:continuous-moments} and \eqref{eq:tag-7.1},
$\|\avg_N\|_{L^2_tH}+\|\avg_N\|_{L^6_tH}$ has uniformly bounded
moments of every finite order. Lemma~\ref{lem:stability} and
H\"older's inequality therefore give
\begin{equation}
 \begin{aligned}
 \|\err_N\|_{L^p(\Omega;C_tH)}
 &\le C\left\|1+\bigl(1+\|\avg_N\|_{L^2_tH}\bigr)
                   \|\avg_N\|_{L^6_tH}^3\right\|_{L^{2p}(\Omega)}
          \|g_N\|_{L^{2p}(\Omega;L^2_tH)}\\
 &\le C_T\|g_N\|_{L^{2p}(\Omega;L^\infty_tH)}
 \le C\lambda_{N+1}^{-\gamma/2}.
 \end{aligned}
 \label{eq:spectral-energy-error}
\end{equation}
Here we used $\|g_N\|_{L^2_tH}\le T^{1/2}\|g_N\|_{L^\infty_tH}$
and \eqref{eq:spectral-flux} with moment order $2p$.

By \eqref{eq:continuous-moments}, \eqref{eq:tag-7.1}, and
$\cH^\beta\hookrightarrow L^4$ for $\beta>1/4$,
$\sup_t\|\avg_N(t)\|_4$ has uniformly bounded moments of every
finite order. Using
$\|\avg_N\err_N\|_{4/3}\le\|\avg_N\|_4\|\err_N\|_2$
and the $L^{4/3}\to L^\infty$ and $L^2\to L^\infty$ gradient
estimates in \eqref{eq:heat-lpq} gives
\begin{align*}
 \|\err_N(t)\|_\infty
 &\le C\int_0^t(t-s)^{-7/8}\|\avg_N(s)\|_4\|\err_N(s)\|_2\dd s
   +C\int_0^t(t-s)^{-3/4}\|g_N(s)\|_2\dd s\\
 &\le C_T\left(\sup_s\|\avg_N(s)\|_4\right)
              \left(\sup_s\|\err_N(s)\|_2\right)
       +C_T\sup_s\|g_N(s)\|_2.
\end{align*}
Taking $L^p(\Omega)$ norms and using H\"older's inequality,
\eqref{eq:spectral-flux}, and \eqref{eq:spectral-energy-error}, we obtain
\begin{align*}
 \|\err_N\|_{L^p(\Omega;C_tL^\infty_x)}
 &\le C_T\|\avg_N\|_{L^{2p}(\Omega;C_tL^4_x)}
          \|\err_N\|_{L^{2p}(\Omega;C_tH)}
       +C_T\|g_N\|_{L^p(\Omega;L^\infty_tH)}\le C\lambda_{N+1}^{-\gamma/2}.
\end{align*}
Since $u-u_N=\err_N+R_N$, \eqref{eq:tag-7.1} gives
\[
 \|u-u_N\|_{L^p(\Omega;C_tL^\infty_x)}
 \le \|\err_N\|_{L^p(\Omega;C_tL^\infty_x)}
     +\|R_N\|_{L^p(\Omega;C_tL^\infty_x)}
 \le C\lambda_{N+1}^{-\gamma/2}.\qedhere
\]
\end{proof}

\section{Moments and convergence of the fully discrete method}\label{sec:tamed-exact}
Throughout this section, $u_0\in\cH^1$ is deterministic.
In this section, we first define the continuous interpolation $\widetilde Y^N$ and
prove uniform moment and time-increment bounds in
Proposition~\ref{prop:tamed-moments}. These bounds and
Lemmas~\ref{lem:linfty-tail} and~\ref{lem:stability} give the
$H$-error estimate in Proposition~\ref{prop:energy-error}.
Finally, the heat-semigroup estimate \eqref{eq:heat-lpq} yields
the $L^\infty$-error bound in
Theorem~\ref{thm:tamed-linfty-grid}.

We use the following estimates to analyze the scheme \eqref{eq:tag-6.0}:
\begin{align}
 \|B_N(v)\|_{\cH^{-1}}&\le C\|v\|_4^2,\label{eq:tag-6.1}\\
 \|\widehat B_{N,\tau}(v)\|_{\cH^{-1}}&\le\tau^{-1/4},\label{eq:tag-6.2}\\
 \|\widehat B_{N,\tau}(v)\|_2
 &\le\lambda_N^{1/2}\tau^{-1/4}\le\sqrt{c_0}\tau^{-1/2},\label{eq:tag-6.3}\\
 \|B_N(v)-\widehat B_{N,\tau}(v)\|_{\cH^{-1}}
 &\le C\tau^{1/4}\|v\|_4^4.\label{eq:tag-6.4}
\end{align}
These estimates follow from integration by parts, \eqref{eq:tag-6.0a},
and the inverse inequality. The estimate \eqref{eq:tag-6.3}
uses the mesh condition \eqref{eq:tag-6.0b}.

Define the continuous interpolation
\begin{equation}
 \widetilde Y^N(t)=S_N(t)P_Nu_0+Z_N(t)
 +\int_0^tS_N(t-s)\widehat B_{N,\tau}
        (\widetilde Y^N(\kappa_\tau(s)))\dd s,                 \label{eq:tag-6.5}
\end{equation}
where $Z_N=P_NZ$, $\kappa_\tau(s)=t_m$ on $[t_m,t_{m+1})$, and
$\kappa_\tau(T)=T$. Splitting the integral at $t_m$ and using the
semigroup property gives, for $0\le r\le\tau$,
\begin{align*}
 \widetilde Y^N(t_m+r)
 &=S_N(r)\widetilde Y^N(t_m)
   +\int_0^r S_N(r-s)\widehat B_{N,\tau}(\widetilde Y^N(t_m))\dd s\\
 &\quad+\int_{t_m}^{t_m+r}S_N(t_m+r-s)\dd W(s).
\end{align*}
At $r=\tau$, this 
recovers \eqref{eq:tag-6.0} and
$\widetilde Y^N(t_m)=Y_m^N$.
Since the constant Fourier coefficient of $\widehat B_{N,\tau}$  is zero, we also have \begin{equation}
 \langle\widetilde Y^N(t),e_0\rangle
 =\langle u_0,e_0\rangle+\langle W(t),e_0\rangle
 =\langle u(t),e_0\rangle.
 \label{eq:discrete-shared-mean}
\end{equation} 

\begin{proposition}
\label{prop:tamed-moments}
For every $1\le p<\infty$ and $0\le\beta<1/2$, uniformly over meshes satisfying
$\tau^{1/2}\lambda_N\le c_0$,
\begin{align}
 \sup_{N,\tau}\E\biggl[\sup_{t\le T}
       \norm{\widetilde Y^N(t)}_{\cH^\beta}^p
 &+\norm{\widetilde Y^N}_{L^8(0,T;L^4)}^p
 \biggr]<\infty.                                      \label{eq:tag-6.6}
\end{align}
Moreover, for every $0<\eta<1/4$,
\begin{equation}
 \left(\E\norm{\widetilde Y^N-
 \widetilde Y^N\circ\kappa_\tau}_{L^4(0,T;L^4)}^p\right)^{1/p}
 \le C\tau^\eta.                                    \label{eq:tag-6.7}
\end{equation}

\end{proposition}

\begin{proof}
We first prove \eqref{eq:tag-6.6} using energy and heat-equation
estimates, and then establish \eqref{eq:tag-6.7} using fractional
Morrey embedding. It suffices to consider $p\ge4$.
For $\alpha\ge1$, set $v^{N,\tau}=\widetilde Y^N-Z_{\alpha,N}$.
On $[t_m,t_{m+1}]$,
\begin{equation}
 \partial_t v^{N,\tau}+\nu Lv^{N,\tau}=\widehat B_{N,\tau}(Y_m^N)+\alpha Z_{\alpha,N},\qquad v^{N,\tau}(0)=P_Nu_0.
 \label{eq:tag-6.8a}
\end{equation}
All constants below are independent of $N,\tau,\alpha$ unless the
$\alpha$ dependence is displayed.

Fix $0<h\le\tau$. The energy identity on $[t_m,t_m+h]$ reads
\begin{align}
 &\|v^{N,\tau}(t_m+h)\|_2^2-\|v^{N,\tau}(t_m)\|_2^2
       +2\nu\int_{t_m}^{t_m+h}\|\partial_xv^{N,\tau}\|_2^2\dd t\notag\\
 &\quad=2h\langle\widehat B_{N,\tau}(Y_m^N),v^{N,\tau}(t_m)\rangle
       +2\alpha\int_{t_m}^{t_m+h}\langle Z_{\alpha,N},v^{N,\tau}\rangle\dd t\notag\\
 &\qquad+2\int_{t_m}^{t_m+h}\langle\widehat B_{N,\tau}(Y_m^N),
                 v^{N,\tau}(t)-v^{N,\tau}(t_m)\rangle\dd t.
 \label{eq:tag-6.8c}
\end{align}
The conservative cancellation
$\int_{\T_\ell}v^{N,\tau}(t_m)^2(v^{N,\tau}(t_m))_x\dd x=0$ and $0<\Theta_\tau\le1$ give
\begin{equation}
 |\langle \widehat B_{N,\tau}(Y_m^N),v^{N,\tau}(t_m)\rangle|
 \le C\|(v^{N,\tau}(t_m))_x\|_2(\|Z_{\alpha,N}(t_m)\|_\infty\|v^{N,\tau}(t_m)\|_2+\|Z_{\alpha,N}(t_m)\|_\infty^2).
 \label{eq:tag-6.8d}
\end{equation}
The mesh restriction \eqref{eq:tag-6.0b} implies
$\tau\mu_N\le c_0$ for $\tau\le1$.
Since $y/(1-e^{-y})\le1+y$ for $y\ge0$ (with value $1$ at zero),
Parseval's identity gives
\begin{equation}
 h\|v_x\|_2^2\le(1+2\nu c_0)\int_0^h\|\partial_xS_N(r)v\|_2^2\dd r,
 \qquad v\in H_N,\quad 0<h\le\tau.\label{eq:tag-6.8e}
\end{equation}
For $0\le r\le h$, the mild formula reads
\begin{align*}
 v^{N,\tau}(t_m+r)
 &=S_N(r)v^{N,\tau}(t_m)
   +\int_0^rS_N(r-s)\widehat B_{N,\tau}(Y_m^N)\dd s\\
 &\quad+\alpha\int_0^rS_N(r-s)Z_{\alpha,N}(t_m+s)\dd s.
\end{align*}
The $L^2$ gradient estimate in \eqref{eq:heat-lpq} bounds the
gradient of the drift integral, while Young's
convolution inequality in time bounds the integral containing
$Z_{\alpha,N}$:
\begin{align*}
 \int_0^h\left\|\int_0^r\partial_xS_N(r-s)
       \widehat B_{N,\tau}(Y_m^N)\dd s\right\|_2^2\dd r
 &\le Ch^2\|\widehat B_{N,\tau}(Y_m^N)\|_2^2,\\
 \int_0^h\left\|\alpha\int_0^r\partial_xS_N(r-s)
       Z_{\alpha,N}(t_m+s)\dd s\right\|_2^2\dd r
 &\le C\alpha^2h\int_{t_m}^{t_m+h}\|Z_{\alpha,N}\|_2^2\dd t.
\end{align*}
Applying \eqref{eq:tag-6.8e} to $v^{N,\tau}(t_m)$ and using
these two bounds gives
\begin{equation}
 \begin{aligned}
 h\|(v^{N,\tau}(t_m))_x\|_2^2
 &\le C\int_{t_m}^{t_m+h}\|\partial_xv^{N,\tau}\|_2^2\dd t
       +Ch^2\|\widehat B_{N,\tau}(Y_m^N)\|_2^2\\
 &\quad+C\alpha^2h\int_{t_m}^{t_m+h}\|Z_{\alpha,N}\|_2^2\dd t.
 \end{aligned}
 \label{eq:tag-6.8f}
\end{equation}
By \eqref{eq:tag-6.3}, the drift term satisfies
$h^2\|\widehat B_{N,\tau}(Y_m^N)\|_2^2\le Ch$.

To estimate the last integral in \eqref{eq:tag-6.8c}, subtract
$v^{N,\tau}(t_m)$ from the mild formula and use
$\|(S_N(r)-I)v\|_2\le Cr^{1/2}\|v_x\|_2$,
\eqref{eq:tag-6.3}, and $h\le\tau$ to obtain
\begin{equation}
 \begin{aligned}
 &\int_{t_m}^{t_m+h}|\langle\widehat B_{N,\tau}(Y_m^N),
             v^{N,\tau}(t)-v^{N,\tau}(t_m)\rangle|\dd t\\
 &\qquad\le Ch\|(v^{N,\tau}(t_m))_x\|_2+Ch
       +C\alpha h\|Z_{\alpha,N}\|_{L^2(t_m,t_m+h;H)}.
 \end{aligned}
 \label{eq:tag-6.8g}
\end{equation}
\begin{samepage}
For every $\varepsilon>0$, \eqref{eq:tag-6.8f} and
\eqref{eq:tag-6.3} bound the grid-gradient term produced by Young's
inequality in \eqref{eq:tag-6.8d} and \eqref{eq:tag-6.8g} as follows:
\[
 \begin{aligned}
 \varepsilon h\|(v^{N,\tau}(t_m))_x\|_2^2
 &\le C\varepsilon\int_{t_m}^{t_m+h}
          \|\partial_xv^{N,\tau}\|_2^2\dd t+C\varepsilon h\\
 &\quad+C\varepsilon\alpha^2h
          \int_{t_m}^{t_m+h}\|Z_{\alpha,N}\|_2^2\dd t.
 \end{aligned}
\]
\end{samepage}
Choosing $\varepsilon$ so that $C\varepsilon\le\nu$ and using
$2\alpha|\langle Z_{\alpha,N},v^{N,\tau}\rangle|
\le\|v^{N,\tau}\|_2^2+\alpha^2\|Z_{\alpha,N}\|_2^2$
in \eqref{eq:tag-6.8c}, we obtain
\begin{align*}
 &\|v^{N,\tau}(t_m+h)\|_2^2+\nu\int_{t_m}^{t_m+h}\|\partial_xv^{N,\tau}\|_2^2\dd t\\
 &\quad\le[1+Ch(1+\|Z_{\alpha,N}(t_m)\|_\infty^2)]\|v^{N,\tau}(t_m)\|_2^2+Ch(1+\|Z_{\alpha,N}(t_m)\|_\infty^4)\\
 &\qquad+C\alpha^2\int_{t_m}^{t_m+h}\|Z_{\alpha,N}\|_2^2\dd t
 +C\int_{t_m}^{t_m+h}\|v^{N,\tau}\|_2^2\dd t.
\end{align*}
By Gronwall's inequality, we obtain
\begin{equation}
 \begin{aligned}
 &\sup_{0\le h\le\tau}\left[
       \|v^{N,\tau}(t_m+h)\|_2^2
       +\nu\int_{t_m}^{t_m+h}\|\partial_xv^{N,\tau}\|_2^2\dd t\right]\\
 &\quad\le[1+C_G\tau(1+\|Z_{\alpha,N}(t_m)\|_\infty^2)]
                    \|v^{N,\tau}(t_m)\|_2^2
       +C\tau(1+\|Z_{\alpha,N}(t_m)\|_\infty^4)+R_m,
 \end{aligned}
 \label{eq:tag-6.8}
\end{equation}
where $C_G\ge1$ is independent of $N,\tau,\alpha$ and
$R_m=C\tau+C(1+\alpha^2)\int_{t_m}^{t_{m+1}}(1+\|Z_{\alpha,N}\|_2^2)\dd t$.
Repeating the above arguments, we obtain 
\begin{align*}
 &\sup_{t\le T}\|v^{N,\tau}(t)\|_2^2
       +\nu\int_0^T\|\partial_xv^{N,\tau}\|_2^2\dd t\\
 &\quad\le C\left(1+\|u_0\|_2^2
       +\tau\sum_{m=0}^{M-1}\|Z_{\alpha,N}(t_m)\|_\infty^4+\sum_{m=0}^{M-1}R_m\right) \exp\!\left(C_GT
       +C_G\tau\sum_{m=0}^{M-1}\|Z_{\alpha,N}(t_m)\|_\infty^2\right).
\end{align*}
Choosing $\alpha\ge\alpha_*(16pC_G)$, we apply
\eqref{eq:tag-4.3}, \eqref{eq:tag-4.1b}, and \eqref{eq:tag-4.5}
together with H\"older's inequality and take square roots to obtain
\begin{equation}
 \left\|\sup_{t\le T}\|v^{N,\tau}(t)\|_2\right\|_{L^{2p}(\Omega)}
 +\|v^{N,\tau}\|_{L^{2p}(\Omega;L^2_t\cH^1_x)}
 \le C_p(1+\alpha).\label{eq:tag-6.9}
\end{equation}

We next estimate the grid values. Taking $h=\tau$ in
\eqref{eq:tag-6.8f} and summing over
$m=0,\ldots,M-1$, we use
\begin{align*}
 \tau\sum_{m=0}^{M-1}\|v^{N,\tau}(t_m)\|_2^2
 &\le T\sup_t\|v^{N,\tau}(t)\|_2^2, \; 
 \sum_{m=0}^{M-1}\tau^2\|\widehat B_{N,\tau}(Y_m^N)\|_2^2
\le CT
\end{align*}
to obtain
\begin{equation}
 \left\|\left(\tau\sum_{m=0}^{M-1}\|v^{N,\tau}(t_m)\|_{\cH^1}^2\right)^{1/2}
 \right\|_{L^{2p}(\Omega)}\le C_p(1+\alpha).
 \label{eq:tag-6.9a}
\end{equation}

We now bound the right-hand side of \eqref{eq:tag-6.8a} in
$L^p(\Omega;L^4_t\cH^{-1}_x)$. By \eqref{eq:GN},
\[
 \tau\sum_{m=0}^{M-1}\|v^{N,\tau}(t_m)\|_4^8
 \le C\sup_t\|v^{N,\tau}(t)\|_2^6\,\tau\sum_{m=0}^{M-1}\|v^{N,\tau}(t_m)\|_{\cH^1}^2.
\]
Since $Y_m^N=v^{N,\tau}(t_m)+Z_{\alpha,N}(t_m)$,
\eqref{eq:tag-6.9}, \eqref{eq:tag-6.9a}, \eqref{eq:tag-4.1b},
and $\cH^{3/8}\hookrightarrow L^4$ give
\begin{equation}
 \left\|\left(\tau\sum_{m=0}^{M-1}\|Y_m^N\|_4^8\right)^{1/8}
 \right\|_{L^{2p}(\Omega)}\le C_p(1+\alpha).
 \label{eq:tag-6.9b}
\end{equation}
For $f(t)=\widehat B_{N,\tau}(Y_m^N)+\alpha Z_{\alpha,N}(t)$ on $[t_m,t_{m+1})$,
\eqref{eq:tag-6.1}, \eqref{eq:tag-6.9b}, and
\eqref{eq:tag-4.5b} yield
\begin{equation}
 \|f\|_{L^p(\Omega;L^4_t\cH^{-1}_x)}
 \le C_p(1+\alpha)^2.\label{eq:tag-6.9c}
\end{equation}
Applying \eqref{eq:tag-4.9} and \eqref{eq:tag-4.11} to
\eqref{eq:tag-6.8a} therefore gives
\begin{equation}
 \|v^{N,\tau}\|_{L^p(\Omega;C_t\cH^\beta_x)}
 +\|v^{N,\tau}\|_{L^p(\Omega;W^{\vartheta,4}_t\cH^{1/4}_x)}
 \le C_{p,\beta,\vartheta}(1+\alpha)^2,
 \quad \beta<1/2,\quad 0<\vartheta<3/8.
 \label{eq:tag-6.9d}
\end{equation}
Also, \eqref{eq:GN} and \eqref{eq:tag-6.9} imply
$\|v^{N,\tau}\|_{L^p(\Omega;L^8_tL^4_x)}\le C_p(1+\alpha)$.
Fixing $\alpha=\alpha_*(16pC_G)$ and adding $Z_{\alpha,N}$ using
\eqref{eq:tag-4.1b} now prove \eqref{eq:tag-6.6}.

Finally, we prove \eqref{eq:tag-6.7}. For $1/4<\vartheta<3/8$,
\eqref{eq:morrey} and $\cH^{1/4}\hookrightarrow L^4$ give
\begin{align*}
 \sum_{m=0}^{M-1}\int_{t_m}^{t_{m+1}}
 \|v^{N,\tau}(t)-v^{N,\tau}(t_m)\|_4^4\dd t
 &\le C\tau^{4\vartheta}\sum_{m=0}^{M-1}
 [v^{N,\tau}]_{W^{\vartheta,4}(t_m,t_{m+1};\cH^{1/4})}^4\\
 &\le C\tau^{4\vartheta}
 [v^{N,\tau}]_{W^{\vartheta,4}(0,T;\cH^{1/4})}^4.
\end{align*}
Using \eqref{eq:tag-6.9d} for $v^{N,\tau}$ and applying Minkowski's
inequality to the Gaussian increment bound \eqref{eq:tag-4.6}, we obtain
\[
 \|\widetilde Y^N-\widetilde Y^N\circ\kappa_\tau\|_{L^p(\Omega;L^4_tL^4_x)}
 \le C_p(\tau^\vartheta+\tau^{1/4})
 \le C_p\tau^\eta,\qquad 0<\eta<1/4.
\]
\end{proof}

\subsection{An intermediate strong error in \texorpdfstring{$H$}{H}}

Set $\widehat Y^N=\widetilde Y^N+R_N$ and $\err=u-\widehat Y^N$.
Then it can be seen that 
\begin{equation}
 \partial_t\err+\nu L\err=\partial_x(\avg\err+g),\qquad
 \avg=\tfrac12(u+\widehat Y^N),\qquad \err(0)=0,
 \label{eq:tag-6.10}
\end{equation}
where $g=g_1+g_2+g_3+g_4$ and
\begin{align}
 g_1&=F(\widehat Y^N)-F(\widetilde Y^N),&
 g_2&=F(\widetilde Y^N)-F(\widetilde Y^N\circ\kappa_\tau),\notag\\
 g_3&=(I-P_N)F(\widetilde Y^N\circ\kappa_\tau),&
 g_4&=(1-\Theta_\tau(\widetilde Y^N\circ\kappa_\tau))P_NF(\widetilde Y^N\circ\kappa_\tau).
 \label{eq:tag-6.11}
\end{align}
For later use choose the exponents
\begin{equation}
 0<\gamma<\delta<\tfrac12,\qquad
 \frac{\delta+1/2}{2}<\beta<\tfrac12,\qquad
 \frac\gamma2<\eta<\tfrac14,\qquad 4<\rho<\frac{8\eta}{\gamma}.\label{eq:tag-7.3b}
\end{equation}

\begin{proposition}\label{prop:energy-error}
Under \eqref{eq:tag-6.0b}, for every finite $p\ge2$,
$0<\gamma<1/2$, and $0<\eta<1/4$,
\begin{equation}
 \|\err\|_{L^p(\Omega;C_tH)}
 +\left\|\max_{0\le m\le M}\|u(t_m)-Y_m^N\|_2\right\|_{L^p(\Omega)}
 \le C(\lambda_{N+1}^{-\gamma/2}+\tau^\eta).
 \label{eq:intermediate-energy}
\end{equation}
\end{proposition}
\begin{proof}
We estimate the four flux errors in \eqref{eq:tag-6.11}.
By \eqref{eq:tag-7.1}, \eqref{eq:tag-6.6}, and H\"older's inequality,
\begin{align*}
 \|g_1\|_{L^p(\Omega;L^2_tH)}
 &\le\|\widetilde Y^N\|_{L^{2p}(\Omega;L^4_{t,x})}
       \|R_N\|_{L^{2p}(\Omega;L^4_{t,x})}\\
 &\quad+\tfrac12\|R_N\|_{L^{2p}(\Omega;L^4_{t,x})}^2
 \le C_p\lambda_{N+1}^{-\gamma/2}.
\end{align*}
For $g_2$, the identity $F(v)-F(w)=(v+w)(v-w)/2$,
H\"older's inequality, and \eqref{eq:tag-6.6}, \eqref{eq:tag-6.9b},
\eqref{eq:tag-6.7} at moment order $2p$ give
\begin{align*}
 \|g_2\|_{L^p(\Omega;L^2_tH)}
 &\le\tfrac12
 \|\widetilde Y^N+\widetilde Y^N\circ\kappa_\tau\|_{L^{2p}(\Omega;L^4_{t,x})}
 \|\widetilde Y^N-\widetilde Y^N\circ\kappa_\tau\|_{L^{2p}(\Omega;L^4_{t,x})}
 \le C_p\tau^\eta.
\end{align*}
Choose $\gamma<\delta<1/2$ and $(\delta+1/2)/2<\beta<1/2$.
Using \eqref{eq:tag-2.3}, \eqref{eq:tag-2.10}, and \eqref{eq:tag-6.6},
\begin{align*}
 \|g_3\|_{L^p(\Omega;L^2_tH)}
 &\le\lambda_{N+1}^{-\delta/2}
       \|F(\widetilde Y^N\circ\kappa_\tau)\|_{L^p(\Omega;L^2_t\cH^\delta)}\\
 &\le C_T\lambda_{N+1}^{-\delta/2}
       \|\widetilde Y^N\|_{L^{2p}(\Omega;C_t\cH^\beta)}^2
 \le C_p\lambda_{N+1}^{-\delta/2}.
\end{align*}
Finally, \eqref{eq:tag-6.0a} and \eqref{eq:tag-6.1} imply
\begin{equation}
 \|(1-\Theta_\tau(v))P_NF(v)\|_2
 \le C\tau^{1/4}\|v\|_4^4.\label{eq:tag-6.11b}
\end{equation}
By \eqref{eq:tag-6.11b} and \eqref{eq:tag-6.9b},
\begin{align*}
 \|g_4\|_{L^p(\Omega;L^2_tH)}
 &\le C\tau^{1/4}
 \left\|\left(\tau\sum_{m=0}^{M-1}\|Y_m^N\|_4^8\right)^{1/8}
 \right\|_{L^{4p}(\Omega)}^4
 \le C_p\tau^{1/4}.
\end{align*}
Combining the estimates for $g_1,\ldots,g_4$ yields
\begin{equation}
 \|g\|_{L^p(\Omega;L^2_tH)}
 \le C(\lambda_{N+1}^{-\gamma/2}+\tau^\eta).
 \label{eq:tag-6.12}
\end{equation}
By \eqref{eq:continuous-moments}, \eqref{eq:tag-7.1}, and
\eqref{eq:tag-6.6}, we have $\err\in C_tH$ and
$\|\avg\|_{L^{8p}(\Omega;C_tH)}\le C_p$.
Testing \eqref{eq:tag-6.10} with $\err$ and using
\eqref{eq:tag-6.11} and Young's inequality gives
\[
 \sup_{t\le T}\|\err(t)\|_2^2
 +\nu\int_0^T\|\partial_x\err\|_2^2\dd t
 \le C\left(\|u\|_{L^4_{t,x}}^4
       +\|\widetilde Y^N\circ\kappa_\tau\|_{L^4_{t,x}}^4\right)<\infty.
\]
Thus $\err\in L^2_t\cH^1$. The zero-mean condition in
Lemma~\ref{lem:stability} follows from \eqref{eq:tag-6.10} and
$\err(0)=0$:
\[
 \int_{\T_\ell}\err(t,x)\dd x
 =\int_{\T_\ell}\err(0,x)\dd x=0.
\]
Applying \eqref{eq:tag-3.2}, H\"older's inequality, and
\eqref{eq:tag-6.12} at moment order $2p$, we obtain
\begin{align*}
 \|\err\|_{L^p(\Omega;C_tH)}
 &\le C\left\|1+(1+\|\avg\|_{L^2_tH})\|\avg\|_{L^6_tH}^3
       \right\|_{L^{2p}(\Omega)}
       \|g\|_{L^{2p}(\Omega;L^2_tH)}\\
 &\le C_p\|g\|_{L^{2p}(\Omega;L^2_tH)}
 \le C_p(\lambda_{N+1}^{-\gamma/2}+\tau^\eta).
\end{align*}
At the grid points $u(t_m)-Y_m^N=\err(t_m)+R_N(t_m)$, so
\eqref{eq:tag-7.1} yields
\begin{align*}
 \left\|\max_{0\le m\le M}\|u(t_m)-Y_m^N\|_2\right\|_{L^p(\Omega)}
 &\le\|\err\|_{L^p(\Omega;C_tH)}
       +\|R_N\|_{L^p(\Omega;C_tH)}\\
 &\le C_p(\lambda_{N+1}^{-\gamma/2}+\tau^\eta).\qedhere
\end{align*}
\end{proof}

\subsection{Proof of Theorem~\ref{thm:tamed-linfty-grid}}
Finally, we prove the $L^\infty$-error bound in
Theorem~\ref{thm:tamed-linfty-grid} by applying the heat-semigroup
estimate \eqref{eq:heat-lpq} to the error equation \eqref{eq:tag-6.10}
and using the $H$-error bound from Proposition~\ref{prop:energy-error}.

\begin{proof}[Proof of Theorem~\ref{thm:tamed-linfty-grid}]
Use $\err,\avg,g_1,\ldots,g_4$ from \eqref{eq:tag-6.10} and
\eqref{eq:tag-6.11}, with the exponents fixed in \eqref{eq:tag-7.3b}.
Since $\widehat Y^N=\widetilde Y^N+R_N$,
\[
 \|g_1(t)\|_2
 \le \|\widetilde Y^N(t)\|_2\|R_N(t)\|_\infty
     +C\|R_N(t)\|_\infty^2.
\]
Proposition~\ref{prop:tamed-moments} and \eqref{eq:tag-7.1} therefore imply
\begin{equation}
 \|g_1\|_{L^p(\Omega;C_tH)}
 \le C\lambda_{N+1}^{-\gamma/2}.                       \label{eq:tag-7.3d}
\end{equation}

At the grid points $Y_m^N=\widetilde Y^N(t_m)$. Estimates \eqref{eq:tag-2.10}, \eqref{eq:tag-2.3},
and \eqref{eq:tag-6.6} give
\[
 \sup_{t\le T}\|g_3(t)\|_2
 \le C\lambda_{N+1}^{-\delta/2}
       \sup_m\|Y_m^N\|_{\cH^\beta}^2,
\]
and hence
\begin{equation}
 \|g_3\|_{L^p(\Omega;L^\infty_tH)}
 \le C\lambda_{N+1}^{-\delta/2}.                           \label{eq:tag-7.3e}
\end{equation}
By $\cH^\beta\hookrightarrow L^4$ and \eqref{eq:tag-6.6}, all finite moments of
$\sup_m\|Y_m^N\|_4$ and $\sup_t\|\widetilde Y^N(t)\|_4$ are uniformly
bounded.  Estimate \eqref{eq:tag-6.11b} thus gives
\begin{equation}
 \|g_4\|_{L^p(\Omega;L^\infty_tH)}
 \le C\tau^{1/4}.                                         \label{eq:tag-7.3f}
\end{equation}

For the freezing error $g_2$, the identity
$F(v)-F(w)=\tfrac12(v+w)(v-w)$ gives
\begin{align*}
 \|g_2\|_{L^4_tH}
 &\le C\left(
 \|\widetilde Y^N\|_{L^\infty_tL^4_x}
 +\|\widetilde Y^N\circ\kappa_\tau\|_{L^\infty_tL^4_x}\right)
 \|\widetilde Y^N-\widetilde Y^N\circ\kappa_\tau\|_{L^4_tL^4_x}.
\end{align*}
Using \eqref{eq:tag-6.6}--\eqref{eq:tag-6.7}, H\"older's inequality, and the choice
$\eta<1/4$, we obtain
\begin{equation}
 \|g_2\|_{L^p(\Omega;L^4_tH)}
 \le C\tau^\eta.                                          \label{eq:tag-7.3g}
\end{equation}
Estimates \eqref{eq:tag-7.3d}--\eqref{eq:tag-7.3f} give convergence
rates in $L^\infty_tH$ for $g_1,g_3,g_4$. The $L^4_tH$ rate
\eqref{eq:tag-7.3g} for $g_2$ cannot be used directly with
\eqref{eq:heat-lpq}, since
\[
 \int_0^T(r^{-3/4})^{4/3}\dd r
 =\int_0^T r^{-1}\dd r=\infty.
\]
To obtain an $L^\rho_tH$ rate with $\rho>4$, we use
\eqref{eq:tag-6.6} and $\cH^\beta\hookrightarrow L^4$ to bound
\begin{equation}
 \sup_{N,\tau}\|g_2\|_{L^p(\Omega;L^\infty_tH)}
 \le C\sup_{N,\tau}\|\widetilde Y^N\|_{L^{2p}(\Omega;C_tL^4_x)}^2
 \le C_p.
 \label{eq:tag-7.3h}
\end{equation}
Bound \eqref{eq:tag-7.3h} alone has no decay in $\tau$.
Interpolating it with \eqref{eq:tag-7.3g} in time and applying
H\"older's inequality gives
\begin{align}
 \|g_2\|_{L^p(\Omega;L^\rho_tH)}
 &\le \|g_2\|_{L^p(\Omega;L^4_tH)}^{4/\rho}
        \|g_2\|_{L^p(\Omega;L^\infty_tH)}^{1-4/\rho}
 \le C\tau^{4\eta/\rho}.\label{eq:freezing-interpolation}
\end{align}

The mild form of \eqref{eq:tag-6.10} is
\[
 \err(t)=\int_0^t\partial_xS(t-s)
 \left[\avg(s)\err(s)+\sum_{i=1}^4g_i(s)\right]\dd s.
\]
We first estimate the integral containing $\avg\err$.
By \eqref{eq:continuous-moments}, \eqref{eq:tag-6.6},
\eqref{eq:tag-7.1}, and $\cH^\beta\hookrightarrow L^4$,
\begin{equation}
 \sup_{N,\tau}\|\avg\|_{L^p(\Omega;C_tL^4_x)}\le C_p
 \qquad\text{for every finite }p.\label{eq:tag-7.3j}
\end{equation}
Since $r^{-7/8}\in L^1(0,T)$, \eqref{eq:heat-lpq}, \eqref{eq:tag-7.3j},
Proposition~\ref{prop:energy-error}, and H\"older's inequality give
\begin{align}
 \left\|\int_0^t
 \partial_xS(t-s)[\avg(s)\err(s)]\dd s\right\|_{L^p(\Omega;C_tL^\infty_x)}
 &\le C_T\|\avg\|_{L^{2p}(\Omega;C_tL^4_x)}
          \|\err\|_{L^{2p}(\Omega;C_tH)}\notag\\
 &\le C_p(\lambda_{N+1}^{-\gamma/2}+\tau^\eta).
 \label{eq:tag-7.3k}
\end{align}

Applying \eqref{eq:heat-lpq} to $g_1,g_3,g_4$ and using
$r^{-3/4}\in L^1(0,T)$, we obtain
\begin{align}
 &\sup_{t\le T}\left\|\int_0^t\partial_xS(t-s)
 [g_1(s)+g_3(s)+g_4(s)]\dd s\right\|_\infty\notag\\
 &\qquad\le C_T\left(
 \|g_1\|_{C_tH}
 +\|g_3\|_{L^\infty_tH}
 +\|g_4\|_{L^\infty_tH}\right).                       \label{eq:tag-7.3l}
\end{align}
For $g_2$, the conjugate exponent $\rho'=\rho/(\rho-1)$ satisfies
$3\rho'/4<1$. Hence
\eqref{eq:heat-lpq} and H\"older's inequality yield
\begin{align}
 \sup_{t\le T}\left\|\int_0^t\partial_xS(t-s)g_2(s)\dd s\right\|_\infty
 &\le C\left(\int_0^T r^{-3\rho'/4}\dd r\right)^{1/\rho'}\|g_2\|_{L^\rho_tH}
 \le C_{T,\rho}\|g_2\|_{L^\rho_tH}.
 \label{eq:tag-7.3m}
\end{align}
Taking $L^p(\Omega)$ norms and using \eqref{eq:freezing-interpolation}
gives a bound $C\tau^{4\eta/\rho}\le C\tau^{\gamma/2}$ by
\eqref{eq:tag-7.3b} and $\tau\le1$.

Combining \eqref{eq:tag-7.3k}, \eqref{eq:tag-7.3l}, and
\eqref{eq:tag-7.3m} with \eqref{eq:tag-7.3d}--\eqref{eq:tag-7.3f} gives
\begin{align*}
 \|\err\|_{L^p(\Omega;C_tL^\infty_x)}
 &\le C\left(\lambda_{N+1}^{-\gamma/2}+\tau^\eta
          +\lambda_{N+1}^{-\delta/2}+\tau^{1/4}+\tau^{4\eta/\rho}\right)\\
 &\le C\left(\lambda_{N+1}^{-\gamma/2}+\tau^{\gamma/2}\right).
\end{align*}
At grid points $u(t_m)-Y_m^N=\err(t_m)+R_N(t_m)$, so
\eqref{eq:tag-7.1} yields \eqref{eq:main-full}:
\begin{align*}
 \left\|\max_{0\le m\le M}\|u(t_m)-Y_m^N\|_\infty\right\|_{L^p(\Omega)}
 &\le\|\err\|_{L^p(\Omega;C_tL^\infty_x)}
       +\|R_N\|_{L^p(\Omega;C_tL^\infty_x)}\\
 &\le C_p(\lambda_{N+1}^{-\gamma/2}+\tau^{\gamma/2}).\qedhere
\end{align*}
\end{proof}

\section{Numerical experiment}\label{sec:experiments}
In this section, we examine the spatial and temporal strong convergence of
scheme~\eqref{eq:tag-6.0}. We take $\ell=2\pi$, $\nu=1$, $T=0.1$, and
$u_0(x)=\sin x$, and use $P=5000$ independent trajectories in each experiment.
Here $Y_m^{N,\tau}$ denotes \eqref{eq:tag-6.0} with step size $\tau$.

At cutoff $N$, both the solution and the noise retain the constant mode
and the first $N$ sine--cosine pairs, giving $2N+1$ real noise modes.
The semigroup and the integral of the frozen drift are evaluated
exactly in each mode. The quadratic term is computed by FFT with
zero-padding to at least $3N+1$ points before projection.
The increment $\zeta_m^{N,\tau}$ from Section~\ref{sec:main} has
independent centered Gaussian coefficients in the orthonormal Fourier basis,
with variance $\tau$ in the constant mode and $(1-e^{-2n^2\tau})/(2n^2)$
in each sine or cosine mode of frequency $1\le n\le N$.
Within each trajectory, all resolutions share the same Brownian coordinates;
coarse increments are obtained by the exact semigroup-weighted sum of fine increments.
We use $\|\cdot\|_{\ell^\infty}$ for the maximum over $2048$ equally
spaced spatial evaluation points in $[0,2\pi)$.

For the spatial experiment, set $M_*=32768$ and $\tau_*=T/M_*$, and compare
$N=4,8,16,32,64$ with the cutoff-$256$ reference, using
\[
 E_N=\left(\frac1P\sum_{i=1}^P\max_{0\le m\le M_*}
 \|Y_m^{N,\tau_*,(i)}-Y_m^{256,\tau_*,(i)}\|_{\ell^\infty}^2\right)^{1/2}.
\]

For the temporal experiment, fix $N=128$ and take $\tau=T/M$,
$M=2048,4096,8192,16384$. The reference solution $\overline Y$ uses
the same exponential update and exact stochastic increments, with the modified taming factor
\[
 \Theta^{\mathrm{ref}}_{\tau_{\mathrm{ref}}}(v)
 =\frac1{1+\tau_{\mathrm{ref}}^{1/2}\|B_{128}(v)\|_{\cH^{-1}}},
 \qquad M_{\mathrm{ref}}=131072,\quad \tau_{\mathrm{ref}}=T/M_{\mathrm{ref}}.
\]
We estimate the $L^2(\Omega)$-error by
\[
 E_\tau=\left(\frac1P\sum_{i=1}^P\max_{0\le m\le M}
 \|Y_m^{128,\tau,(i)}-
 \overline Y_{mM_{\mathrm{ref}}/M}^{(i)}\|_{\ell^\infty}^2\right)^{1/2}.
\]
Both solutions use the same Galerkin space and are compared at common
times, so the temporal experiment concerns the Galerkin equation at $N=128$.
All tested and reference meshes satisfy \eqref{eq:tag-6.0b} with $c_0=128$.

\begin{figure}[htbp]
 \centering
 \includegraphics[width=\linewidth]{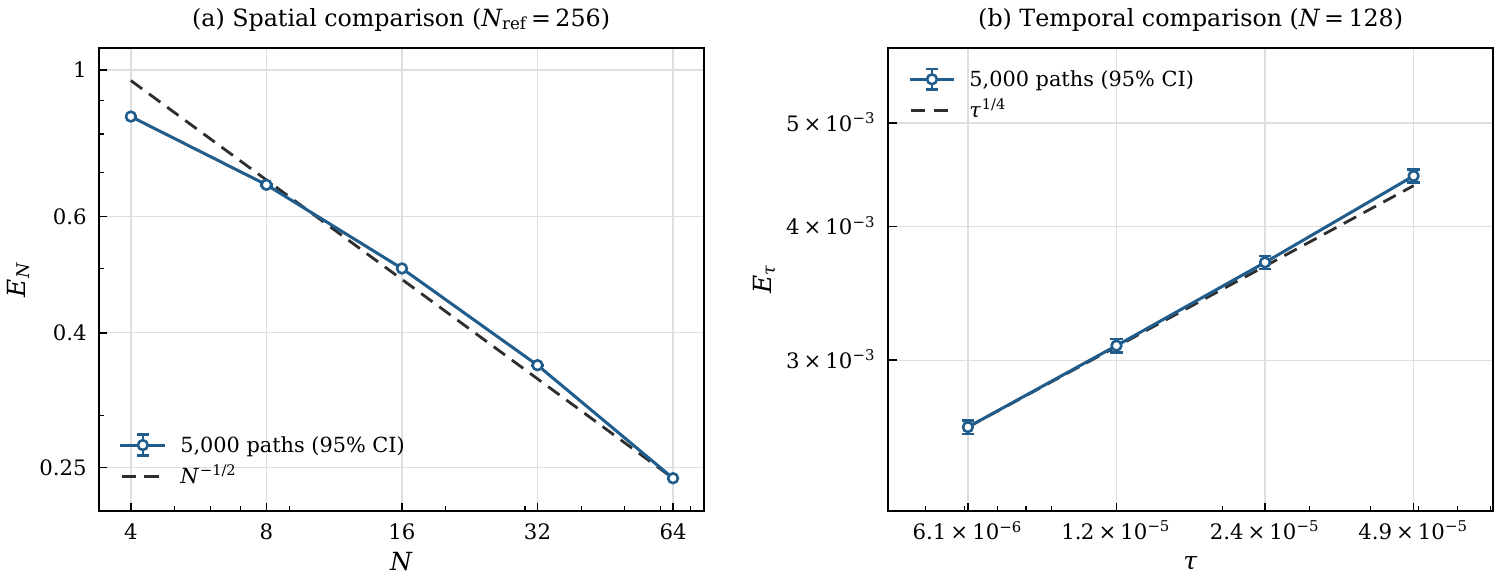}
 \caption{Spatial (left) and temporal (right) $L^2(\Omega)$-errors from
 $P=5000$ trajectories in each experiment. Error bars show pointwise 95\%
 bootstrap intervals. Dashed lines show $N^{-1/2}$ and $\tau^{1/4}$,
 normalized at the finest displayed level.}
 \label{fig:strong}
\end{figure}

Figure~\ref{fig:strong} shows observed strong convergence rates of approximately
$0.455$ in space and $0.260$ in time, close to $1/2$ and $1/4$,
respectively. These numerical results are consistent with the theoretical
estimates in Theorem~\ref{thm:tamed-linfty-grid}.
\FloatBarrier

\section*{Acknowledgement}
The first author thanks Arnulf Jentzen for helpful discussions on this problem. 
The research of the first author is partially supported by the Hong Kong
Research Grants Council GRF grant 15301025, NSFC grant 12522119,
NSFC/RGC Joint Research Scheme N$\_$PolyU5141/24, MOST National Key
$R\&D$ Program No.\ 2024YFA1015900, the internal grants
(P0045336, P0046811) from Hong Kong Polytechnic University and the
CAS AMSS-PolyU Joint Laboratory of Applied Mathematics.
This work is also supported by the National Natural Science Foundation of China
(No.~12471386, No.~12461160278). The authors used GPT to polish the English writing and assist with the
numerical experiments, including checks of the numerical implementation.
The authors assume responsibility for all content.
\par

\begingroup
\small
\bibliographystyle{abbrv}
\bibliography{references}
\endgroup
\end{document}